\documentclass[11pt]{amsart}

\author{Yakov Berchenko-Kogan}
\title{Duality in finite element exterior calculus}

\usepackage{amsmath,amssymb,amsthm,tikz,hyperref}
\usetikzlibrary{cd}

\newtheorem{thm}{Theorem}[section]
\newtheorem{prop}[thm]{Proposition}
\newtheorem{lemma}[thm]{Lemma}
\newtheorem{cor}[thm]{Corollary}

\theoremstyle{definition}
\newtheorem{dfn}[thm]{Definition}

\theoremstyle{remark}
\newtheorem{eg}[thm]{Example}
\newtheorem{question}[thm]{Question}

\newcommand{\mc}{\mathcal}
\newcommand{\mbf}{\mathbf}
\newcommand{\mf}{\mathfrak}
\newcommand{\mb}{\mathbb}
\newcommand{\pp}[2][]{\frac{\partial{#1}}{\partial{#2}}}
\newcommand{\cl}{\overline}
\newcommand{\cO}{\cl{\mbf O}}
\newcommand{\con}{\cO_{\neq0}}
\newcommand{\mr}{\mathring}
\newcommand{\vol}{\mathrm{vol}}
\newcommand{\mg}{\mathbf\Gamma}

\DeclareMathOperator{\image}{image}
\DeclareMathOperator{\Span}{span}

\begin{document}

\begin{abstract}
  In order to generalize finite element methods to differential forms, Arnold, Falk, and Winther constructed two families of spaces of polynomial differential forms on a simplex $T$, the $\mc P_r\Lambda^k(T)$ spaces and the $\mc P_r^-\Lambda^k(T)$ spaces, where $k$ is the degree of the form and $r$ is the degree of its coefficients. 
  The geometric decomposition for these finite element spaces hinges on a duality relationship between the $\mc P$ and $\mc P^-$ spaces proved by Arnold, Falk, and Winther. In this article, we give a natural alternate construction of the $\mc P_r\Lambda^k(T)$ and $\mc P_r^-\Lambda^k(T)$ spaces, leading to a new basis-free proof of this duality relationship using a modified Hodge star operator.
\end{abstract}

\maketitle

\section{Introduction}
Finite element methods are a tool for finding approximate solutions to partial differential equations by triangulating the domain into \emph{elements} and then finding an approximate solution that is a polynomial of degree at most $r$ on each element. Requiring interelement continuity of the solution imposes constraints on these polynomials. By associating each constraint to a shared vertex, edge, or face, one obtains a \emph{geometric decomposition} for the finite element space. That is, the geometric decomposition associates some number of degrees of freedom to each vertex, edge, etc., of an element. Assigning a real number to each of these degrees of freedom uniquely determines the polynomial function on that element, and, if two elements intersect at a subsimplex, the interelement continuity of the polynomial function is equivalent to assigning the same numbers to the degrees of freedom on that subsimplex.

Finite element exterior calculus is the extension of these methods to differential forms \cite{afw06,afw10}. Generalizing finite element methods to differential forms has applications to Maxwell's equations \cite{h02}, elasticity, the Hodge Laplacian, and other problems.
In \cite{afw06,afw10}, Arnold, Falk, and Winther construct the $\mc P_r\Lambda^k(T)$ and $\mc P_r^-\Lambda^k(T)$ spaces of differential forms on a simplex $T$ with polynomial coefficients. They provide a geometric decomposition for these spaces that relies on a duality relationship between the $\mc P$ and $\mc P^-$ spaces. In this article, we give an alternate construction of these families of spaces, leading to a new basis-free proof of this duality relationship.

As discussed in \cite{afw06,afw10}, it is convenient to use barycentric coordinates on $T$. In other words, denote the first orthant by $\cO=\{x\in\mb R^{n+1}\mid x_i\ge0\}$, and let $T=\{x\in\cO\mid x_1+\dotsb+x_{n+1}=1\}$. In Section \ref{prelimsec}, we define a splitting of $\Lambda^k(\cO)$ into vertical and horizontal $k$-forms. A $k$-form $\alpha$ is vertical if $\alpha(X_1,\dotsc,X_k)=0$ whenever the $X_i$ are all parallel to $T$. A $k$-form is horizontal if it is orthogonal to the space of vertical $k$-forms with respect to a nonstandard inner product $g$ that we define.


In Section \ref{prprmsec}, we define $\mbf P_r\Lambda^k(\cO)$ as the space of vertical $k$-forms whose coefficients are homogeneous polynomials of degree $r$, and $\mbf P_r^-\Lambda^k(\cO)$ as the space of horizontal $k$-forms whose coefficients are homogeneous polynomials of degree $r$. We show that the $\mbf P_r\Lambda^k(\cO)$ and $\mbf P_r^-\Lambda^k(\cO)$ spaces are isomorphic via appropriate restriction and contraction maps to the $\mc P_r\Lambda^k(T)$ and $\mc P_r^-\Lambda^k(T)$ spaces of Arnold, Falk, and Winther \cite{afw06,afw10}.

Via the exterior derivative, the $\mc P_r$ and $\mc P_r^-$ spaces each form a cochain complex. We discuss the corresponding cochain complex for the $\mbf P_r$ and $\mbf P_r^-$ spaces in Section \ref{diffsec}.

The heart of this paper is Section \ref{integrationsec}. In this section, we develop a modified Hodge star operator and then use this operator to obtain a quick proof in Theorem \ref{pairingthm} of the duality relationship between the $\mbf P$ and $\mbf P^-$ spaces. Via the correspondence in Section \ref{prprmsec}, we immediately obtain as Corollary \ref{afwcor} the duality relationship between the $\mc P$ and $\mc P^-$ spaces established in \cite{afw10}. 


Finally, in Section \ref{futuresec}, we briefly mention two potential future directions for this work.

\section{Preliminaries}\label{prelimsec}
\begin{dfn}
  Let $\mbf O$ be the first orthant, that is, the set of points in $\mathbb R^{n+1}$ all of whose coordinates are strictly positive. Let $\cO$ denote the set of points in $\mathbb R^{n+1}$ all of whose coordinates are non-negative. Let $\con$ denote $\cO$ excluding the origin.
\end{dfn}
\begin{dfn}
Let
  \begin{align*}
    s&=x_1+x_2+\dotsb+x_{n+1}.
  \end{align*}
  Let $T$ be the $n$-dimensional simplex in $\cO$ defined by the equation $s=1$. 
\end{dfn}

\begin{dfn}
  Let $g$ denote a nonstandard inner product on $T^*_x\cO$ defined by
  \begin{align*}
    \langle dx_i,dx_i\rangle_g&=x_i,&\langle dx_i,dx_j\rangle_g=0\text{ when }i\neq j.
  \end{align*}
\end{dfn}
The inner product $g$ is degenerate when $x$ is on the boundary of $\cO$ but is nondegenerate when $x$ is in $\mbf O$. By extending $g$ to tensor powers of $T^*_x\cO$, we obtain an inner product on $k$-forms $\Lambda^k(\cO)$. On decomposable $k$-forms $\alpha=\alpha_1\wedge\dotsb\wedge\alpha_k$ and $\beta=\beta_1\wedge\dotsb\wedge\beta_k$, this inner product can be computed explicitly by
\begin{equation}\label{innerproductdet}
  \langle\alpha,\beta\rangle_g=\det\bigl(\langle\alpha_i,\beta_j\rangle_g\bigr),
\end{equation}
where the $\alpha_i$ and $\beta_j$ are in $T^*_x\cO$.

\begin{eg}
  \begin{equation*}
    \langle x\,dy\wedge dz,y\,dy\wedge dz\rangle_g=xy\,\langle dy\wedge dz,dy\wedge dz\rangle_g=(xy)(yz)=xy^2z.
  \end{equation*}
\end{eg}

Naturally, $g$ also induces an inner product on $T_x\mbf O$ via $\left\langle \pp{x_i},\pp{x_i}\right\rangle=\frac1{x_i}$ and $\left\langle\pp{x_i},\pp{x_j}\right\rangle=0$ for $i\neq j$.

\subsection{The operators $(ds\wedge)$ and $i_X$}

\begin{dfn}
  Let $X$ be the vector field dual to the one-form $ds$ with respect to $g$. That is, $X$ is defined by the equation
  \begin{equation*}
    \alpha(X)=\langle\alpha,ds\rangle_g\text{ for all }\alpha\in\Lambda^1(\cO).
  \end{equation*}
In other words, $X$ is the gradient of $s$ with respect to the metric $g$. For the gradient of $s$ with respect to the standard Euclidean metric, we will use the standard notation $\nabla s$. Note that $X$ is well-defined even on the boundary of $\cO$ where $g$ is degenerate.

\end{dfn}

\begin{prop}
  The vector field $X$ dual to $ds$ with respect to $g$ is equal to
  \begin{equation*}
    X=\sum_{i=1}^{n+1}x_i\pp{x_i}.
  \end{equation*}
  \begin{proof}
    Observe that $ds=\sum_{i=1}^{n+1}dx_i$, so
    \begin{equation*}
      \langle dx_i,ds\rangle_g=\langle dx_i,dx_i\rangle_g=x_i=dx_i(X).
    \end{equation*}
    Extending by linearity, we see that if $\alpha=\sum_{i=1}^{n+1}\alpha_i\,dx_i$ for scalar functions $\alpha_i$, then
    \begin{equation*}
      \langle\alpha,ds\rangle_g=\alpha(X),
    \end{equation*}
    as desired.
  \end{proof}

\end{prop}

\begin{dfn}
  Let
  \begin{equation*}
    (ds\wedge)\colon \Lambda^k(\cO)\to\Lambda^{k+1}(\cO)
  \end{equation*}
  denote the map
  \begin{equation*}
    (ds\wedge)\colon\alpha\mapsto ds\wedge\alpha.
  \end{equation*}
\end{dfn}

\begin{dfn}
  Let
  \begin{equation*}
    i_X\colon\Lambda^{k+1}(\cO)\to\Lambda^k(\cO)
  \end{equation*}
  be the contraction of $(k+1)$-forms with the vector field $X$.
\end{dfn}
Note that both of these operations are tensorial, that is, they are pointwise operations. We will also need the corresponding right wedge and right contraction operations.
\begin{dfn}
  Let
  \begin{align*}
    (\wedge ds)&\colon\Lambda^k(\cO)\to\Lambda^{k+1}(\cO),&j_X&\colon\Lambda^{k+1}(\cO)\to\Lambda^k(\cO)
  \end{align*}
  denote $(-1)^k(ds\wedge)$ and $(-1)^ki_X$, respectively.
\end{dfn}
As the notation suggests, $(\wedge ds)\colon\alpha\mapsto\alpha\wedge ds$, and $j_X\beta$ inserts $X$ into the \emph{rightmost} slot of $\beta$.

Because the one-form $ds$ is dual to the vector field $X$ with respect to the metric $g$, we have the following standard algebraic fact.
\begin{prop}\label{dsixadjoint}
  The operators $(ds\wedge)$ and $i_X$ are adjoints with respect to $g$. That is, for $\alpha\in\Lambda^k(\cO)$ and $\beta\in\Lambda^{k+1}(\cO)$, we have, at every $x\in\cO$,
  \begin{equation*}
    \langle ds\wedge\alpha,\beta\rangle_g=\langle\alpha,i_X\beta\rangle_g.
  \end{equation*}
  \begin{proof}
    Because of linearity, it suffices to prove the statement for $\beta=\beta_0\wedge\dotsb\wedge\beta_k$, where the $\beta_i$ are one-forms. Let
    \begin{equation*}
      \hat\beta_i=\beta_0\wedge\dotsb\wedge\beta_{i-1}\wedge\beta_{i+1}\wedge\dotsb\wedge\beta_k.
    \end{equation*}
    Then, using the formula for the inner product in equation \eqref{innerproductdet} and expanding the determinant by minors, we find that
    \begin{equation*}
      \begin{split}
      \langle ds\wedge\alpha,\beta\rangle_g&=\sum_{i=0}^{k}(-1)^{i}\langle ds,\beta_i\rangle_g\langle\alpha,\hat\beta_i\rangle_g\\
      &=\sum_{i=0}^k(-1)^i\beta_i(X)\langle\alpha,\hat\beta_i\rangle_g\\
      &=\sum_{i=0}^k\langle\alpha,i_X\beta\rangle_g.\qedhere
      \end{split}
    \end{equation*}
  \end{proof}
\end{prop}
Note that the case where $k=0$ is simply the definition of $X$.

\begin{prop}\label{dsixs}
  We have
  \begin{equation*}
    (ds\wedge)\circ i_X+i_X\circ(ds\wedge) = s.
  \end{equation*}
  That is, for any $\alpha\in\Lambda^k(\cO)$, we have
  \begin{equation}\label{dsixseq}
    ds\wedge i_X\alpha+i_X(ds\wedge\alpha)=s\alpha.
  \end{equation}
  Similarly,
  \begin{equation*}
    (\wedge ds)\circ j_X+j_X\circ(\wedge ds) = s.
  \end{equation*}
  \begin{proof}
    Contraction satisfies a signed Leibniz rule, so we have that
    \begin{equation*}
      i_X(ds\wedge\alpha)=(i_Xds)\alpha-ds\wedge i_X\alpha.
    \end{equation*}
    Meanwhile, since $X=\sum_{i=1}^{n+1}x_i\pp{x_i}$ and $ds=\sum_{i=1}^{n+1}dx_i$, we see that $i_Xds=\sum_{i=1}^{n+1}x_i=s$. The first claim follows. The second claim follows from the facts that
    \begin{align*}
      ((\wedge ds)\circ j_X)\alpha&=((-1)^{k-1}(ds\wedge)\circ (-1)^{k-1}i_X)\alpha=((ds\wedge)\circ i_X)\alpha,\\
      (j_X\circ(\wedge ds))\alpha&=((-1)^ki_X\circ(-1)^k(ds\wedge))\alpha=(i_X\circ(ds\wedge))\alpha.\qedhere
    \end{align*}
  \end{proof}
\end{prop}

\begin{cor}\label{dsixexact}
  On forms on $\con$, the following sequences are exact.
  \begin{equation*}
    \begin{tikzcd}[row sep=0]
      \Lambda^{k-1}(\con)\arrow[r,"ds\wedge"]&\Lambda^k(\con)\arrow[r,"ds\wedge"]&\Lambda^{k+1}(\con)\\
      \Lambda^{k-1}(\con)&\Lambda^k(\con)\arrow[l,"i_X" above]&\Lambda^{k+1}(\con)\arrow[l,"i_X" above]
    \end{tikzcd}
  \end{equation*}
  In other words, $\ker(ds\wedge)=\image(ds\wedge)$ and $\ker i_X=\image i_X$.
  \begin{proof}
    Antisymmetry implies that $ds\wedge ds\wedge\alpha=0$ and $i_Xi_X\alpha=0$, so it remains to show that if $ds\wedge\alpha=0$, then $\alpha=ds\wedge\beta$ for some $\beta$, and if $i_X\alpha=0$, then $\alpha=i_X\beta$ for some $\beta$. These are standard algebraic facts, but it will be useful to show how they follow from Proposition \ref{dsixs} using the nonvanishing of $s$ on $\con$.

    Dividing equation \eqref{dsixseq} by $s$, we have
    \begin{equation*}
      s^{-1}ds\wedge i_X\alpha+s^{-1}i_X(ds\wedge\alpha)=\alpha.
    \end{equation*}
    Thus, if $ds\wedge\alpha=0$, then
    \begin{equation}\label{dsixexacteq1}
      \alpha=s^{-1}ds\wedge i_X\alpha=ds\wedge(s^{-1}i_X\alpha).
    \end{equation}
    Likewise, if $i_X\alpha=0$, then
    \begin{equation}\label{dsixexacteq2}
      \alpha=s^{-1}i_X(ds\wedge\alpha)=i_X(s^{-1}ds\wedge \alpha).\qedhere
    \end{equation}
  \end{proof}
\end{cor}

Note that, in fact, the first sequence is exact on all of $\cO$ including the origin, but the second sequence is not, due to the vanishing of $X$ at the origin.


\subsection{Vertical and horizontal spaces of forms}

Thinking of the simplices $\{x\in\con\mid s=\text{const.}\}$ as ``horizontal,'' we can think of $\ker(ds\wedge)=\image(ds\wedge)$ and $\ker i_X=\image i_X$ as \emph{vertical} and \emph{horizontal} subspaces of forms in $\Lambda^k(\con)$, respectively, because the vertical subspace $\image(ds\wedge)$ vanishes when restricted to one of these simplices, and as we will see the spaces of vertical forms and horizontal forms are complements. Bearing in mind Corollary \ref{dsixexact}, we thus make the following choice of notation.

\begin{dfn}
  Let
  \begin{align*}
    \Lambda^k(\con)^\perp&=ds\wedge\Lambda^{k-1}(\con)=\ker(ds\wedge),\\
    \Lambda^k(\con)^\top&=i_X\Lambda^{k+1}(\con)=\ker i_X.
  \end{align*}
\end{dfn}

Note that these subspaces are defined pointwise. The following proposition shows that the vertical and horizontal subspaces are complements in $\Lambda^k(\con)$, using Proposition \ref{dsixs} and the invertibility of $s$ on $\con$.
\begin{prop}\label{dsixsplitting}
  We have the following splitting of $\Lambda^k(\con)$.
  \begin{equation*}
    \Lambda^k(\con)=\Lambda^k(\con)^\perp\oplus\Lambda^k(\con)^\top.
  \end{equation*}
  This splitting is pointwise orthogonal with respect to the inner product $g$.
  \begin{proof}
    We first show orthogonality. By Proposition \ref{dsixadjoint}, the operators $(ds\wedge)$ and $i_X$ are adjoints with respect to the inner product $g$. Therefore,
    \begin{equation*}
      \left\langle ds\wedge\beta,i_X\gamma\right\rangle_g=\left\langle ds\wedge ds\wedge\beta,\gamma\right\rangle_g=0
    \end{equation*}
    for any $ds\wedge\beta\in\Lambda^k(\con)^\perp$ and $i_X\gamma\in\Lambda^k(\con)^\top$. However, because the inner product $g$ is degenerate on the boundary of $\cO$, we cannot immediately conclude that these two spaces have zero intersection.

    Since $s$ is nonzero on $\con$, Proposition \ref{dsixs} tells us
    that
    \begin{equation*}
      s^{-1}\left(ds\wedge i_X\alpha+i_X(ds\wedge\alpha)\right)=\alpha=ds\wedge i_X(s^{-1}\alpha)+i_X(ds\wedge s^{-1}\alpha).
    \end{equation*}
    The left equation tells us that if $\alpha\in\Lambda^k(\con)^\perp=\ker(ds\wedge)$ and $\alpha\in\Lambda^k(\con)^\top=\ker i_X$, then $\alpha=0$. The right equation tells us that any $\alpha$ can be expressed as a sum of a form in $\Lambda^k(\con)^\perp$ and a form in $\Lambda^k(\con)^\top$, as desired.
  \end{proof}
\end{prop}

\begin{prop}\label{vhiso}
  The vertical and horizontal forms are pointwise isomorphic via either of the following two pairs of inverse maps.
    \begin{equation*}
      \begin{tikzcd}
        \Lambda^{k+1}(\con)^\perp \arrow[r, shift left, "i_X"] &\Lambda^k(\con)^\top \arrow[l, shift left, "s^{-1}ds\wedge"]&
        \Lambda^{k+1}(\con)^\perp \arrow[r, shift left, "s^{-1}i_X"] &\Lambda^k(\con)^\top \arrow[l, shift left, "ds\wedge"]
      \end{tikzcd}
    \end{equation*}

  \begin{proof}
    See equations \eqref{dsixexacteq1} and \eqref{dsixexacteq2}.
  \end{proof}
\end{prop}
For reasons of compability with the exterior derivative that will be discussed in Section \ref{diffsec}, we will also need these isomorphisms in terms of the \emph{right} wedge operation $(\wedge ds)=(-1)^k(ds\wedge)$ on $k$-forms and its adjoint, the right contraction $j_X=(-1)^ki_X$ on $(k+1)$-forms.
\begin{cor}
  The vertical and horizontal forms are pointwise isomorphic via either of the following two pairs of inverse maps.
  \begin{equation*}
    \begin{tikzcd}
      \Lambda^{k+1}(\con)^\perp \arrow[r, shift left, "j_X"] &\Lambda^k(\con)^\top \arrow[l, shift left, "s^{-1}(\wedge ds)"]&
      \Lambda^{k+1}(\con)^\perp \arrow[r, shift left, "s^{-1}j_X"] &\Lambda^k(\con)^\top \arrow[l, shift left, "\wedge ds"]
    \end{tikzcd}
  \end{equation*}
\end{cor}

\subsection{Differential forms on the simplex $T$}
We will show that, at each point $x\in T$, the horizontal forms on $\con$ are isomorphic via restriction to forms on the simplex $T$.
\begin{prop}\label{htiso}
  Let $T$ denote the standard $n$-simplex $\{x\in\cO\mid s=1\}$. Let $i\colon T\hookrightarrow\con$ denote the inclusion, and let $i^*\colon\Lambda^k(\con)\to\Lambda^k(T)$ denote the restriction of forms. Let $x\in T$. The restriction map
  \begin{equation*}
    i^*\colon\Lambda^k_x(\con)\to\Lambda^k_x(T)
  \end{equation*}
  is surjective with kernel $\Lambda^k_x(\con)^\perp$.
  \begin{proof}
    Let $\alpha\in\Lambda^k_x(\con)$. If $\alpha\in\Lambda^k_x(\con)^\perp$, then $\alpha=ds\wedge\beta$ for some $\beta$. Then $i^*\alpha=i^*ds\wedge i^*\beta=0$, because $ds(Y)=0$ for any vector $Y$ tangent to $T$. Thus, the kernel of $i^*$ contains $\Lambda^k_x(\con)^\perp$.

    Next we show that the kernel of $i^*$ is not larger than $\Lambda^k_x(\con)^\perp$. Because $\Lambda^k_x(\con)^\perp$ and $\Lambda^k_x(\con)^\top$ are complements, it suffices to show that if $\alpha\in\Lambda^k_x(\con)^\top$ and $i^*\alpha=0$, then $\alpha=0$. Note that $X_x$ is not tangent to $T$. Indeed, $ds(X_x)=s_x=1\neq0$. Consider arbitrary vectors $X_1,\dotsc,X_k$. Because the tangent space of $T$ has codimension one in $\mb R^{n+1}$ and the vector $X_x$ is not tangent to $T$, we know that $X_i=Y_i+c_iX_x$, where $Y_i$ is tangent to $T$ and $c_i$ is a scalar. One can compute $c_i$ explicitly using $ds(X_i)=ds(Y_i)+c_ids(X_x)=c_i$. We thus have that
    \begin{equation*}
      \alpha(X_1,\dotsc,X_k)=\alpha(Y_1+c_1X_x,\dotsc,Y_k+c_kX_x).
    \end{equation*}
    Using multilinearity, we could expand this expression into $2^k$ terms. One term, $\alpha(Y_1,\dotsc,Y_k)$, is zero by the assumption that $i^*\alpha=0$, since $Y_1,\dotsc,Y_k$ are tangent to $T$. The remaining terms all involve plugging in $X_x$ into at least one slot of $\alpha$, and hence they vanish by the assumption that $\alpha\in\Lambda^k_x(\con)^\top$, meaning $i_{X_x}\alpha=0$. We conclude that $\alpha(X_1,\dotsc,X_k)=0$ for arbitrary vectors $X_1,\dotsc,X_k$, and hence $\alpha=0$, as desired.

    Finally, we show the algebraic fact that the restriction map is surjective on the exterior algebra. Given $a\in\Lambda^k_x(T)$, define $\alpha\in\Lambda^k_x(\con)^\top$ by the equation
    \begin{equation*}
      \alpha(X_1,\dotsc,X_k)=a(Y_1,\dotsc,Y_k),
    \end{equation*}
    where the $Y_i$ are the projections of $X_i$ to the tangent space of $T$ given by the equation $Y_i=X_i-c_iX_x=X_i-\left(ds(X_i)\right)X_x$.

    Let $\pi_x\colon\mathbb R^{n+1}\to T_xT$ denote this projection, that is,
    \begin{equation*}
      \pi_x(Z)=Z-\left(ds(Z)\right)X_x.
    \end{equation*}
    Note that if $Y$ is tangent to $T$, then $ds(Y)=0$, so, indeed, $\pi_x(Y)=Y$.

    It is clear that $\alpha$ is multlinear and antisymmetric, and it restricts to $a$. Indeed,
    \begin{equation*}
      i^*\alpha(Y_1,\dotsc,Y_k)=\alpha(Y_1,\dotsc,Y_k)=a(\pi_x(Y_1),\dotsc,\pi_x(Y_k))=a(Y_1,\dotsc,Y_k),
    \end{equation*}
    as desired.
  \end{proof}
\end{prop}
Note that the extension $\alpha$ we constructed in the above proof satisfies $i_{X_x}\alpha=0$, because $\pi_x(X_x)=X_x-(ds(X_x))X_x=X_x-X_x=0$. We let $\pi_x^*$ denote this extension map sending $a$ to $\alpha$.
\begin{dfn}
  For $x\in T$, let $\pi_x^*\colon\Lambda_x^k(T)\to\Lambda_x^k(\con)^\top$ be defined by
  \begin{equation*}
    (\pi_x^*a)(X_1,\dotsc,X_k)=a(\pi_xX_1,\dotsc,\pi_xX_k),
  \end{equation*}
  where $X_1,\dotsc,X_k\in\mb R^{n+1}$ and $\pi_x\colon\mb R^{n+1}\to T_xT$ denotes the map
  \begin{equation*}
    \pi_x(Z)=Z-\left(ds(Z)\right)X_x.
  \end{equation*}
\end{dfn}
Note that $X_x$ is orthogonal to the tangent space of $T$ with respect to the nonstandard inner product $g_x$ on $T_x\mbf O$ via duality: $\langle X_x,Y\rangle_g=ds(Y)=0$. Thus, $\pi_x\colon\mb R^{n+1}\to T_xX$ is, in fact, the $g$-orthogonal projection.

We can thus state the following corollary of Propositions \ref{htiso}:

\begin{cor}\label{htcor}
  At each point $x\in T$, the space of horizontal forms on $\con$ is isomorphic to the space of forms on $T$ via the restriction map $i^*$ and its inverse $\pi^*_x$.
  \begin{equation*}
    \begin{tikzcd}
      \Lambda^k_x(\con)^\top \arrow[r, shift left, "i^*"] &\Lambda^k_x(T). \arrow[l, shift left, "\pi^*_x"] 
    \end{tikzcd}
  \end{equation*}
\end{cor}

Arnold, Falk, and Winther \cite{afw06,afw10} reference the Koszul operator $\kappa\colon\Lambda_x^{k+1}(T)\to\Lambda_x^k(T)$ \cite{gs99}. We will need to understand the Koszul operator in this context, and we will need to understand map $\Lambda^{k+1}_x(\con)^\top\to\Lambda^k_x(\con)^\top$ that corresponds to $\kappa$ via the above isomorphism.

\begin{dfn}\label{koszulvfdfn}
  On $\cO$, let the vector field $X^\kappa$ be defined by
  \begin{equation*}
    X^\kappa=X-\tfrac1{n+1}s\nabla s.
  \end{equation*}
\end{dfn}
The vector field $X^\kappa$ is the orthogonal projection with respect to the \emph{standard} inner product of $X$ to the tangent space of $T$. Indeed, $\nabla s$ is orthogonal to the tangent space of $T$, and
\begin{equation*}
  ds(X^\kappa)=ds(X)-\tfrac1{n+1}s\,ds(\nabla s)=s-\tfrac1{n+1}s(n+1)=0.
\end{equation*}
As such, $X^\kappa$ is the radial tangent vector field on $T$ representing the displacement from the center of the simplex $T$, which is precisely the vector field defining the Koszul operator discussed in \cite{afw06,afw10} with origin at the center of the simplex.

\begin{dfn}\label{koszuldfn}
  Define the operator $\kappa$ to be contraction with the vector field $X^\kappa$. That is, let $\kappa=i_{X^\kappa}$.
\end{dfn}
For $x\in T$, because $X^\kappa$ is tangent to $T$, the operator $\kappa$ is defined in both of the following contexts:
\begin{align*}
  \kappa_x\colon\Lambda^{k+1}_x(\con)&\to\Lambda^k_x(\con),&\kappa_x\colon\Lambda^{k+1}_x(T)&\to\Lambda^k_x(T).
\end{align*}

\begin{prop}\label{koszulcommutes}
  The following diagram commutes
  \begin{equation*}
    \begin{tikzcd}
      \Lambda^{k+1}_x(\con)^\top \arrow[r, shift left, "i^*"] \arrow[d, "\kappa_x"]&\Lambda^{k+1}_x(T). \arrow[l, shift left, "\pi^*_x"] \arrow[d,"\kappa_x"] \\
      \Lambda^k_x(\con)^\top \arrow[r, shift left, "i^*"] &\Lambda^k_x(T). \arrow[l, shift left, "\pi^*_x"] 
    \end{tikzcd}
  \end{equation*}
  \begin{proof}
    Let $\alpha\in\Lambda^{k+1}_x(\con)^\perp$. We follow definitions, checking for tangent vectors $Y_1,\dotsc,Y_k$ that
    \begin{multline*}
      \kappa_x(i^*\alpha)(Y_1,\dotsc,Y_k)=i^*\alpha(X^\kappa_x,Y_1,\dotsc,Y_k)=\alpha(X^\kappa_x,Y_1,\dotsc,Y_k)\\
      =\kappa_x\alpha(Y_1,\dotsc,Y_k)=i^*(\kappa_x\alpha)(Y_1,\dotsc,Y_k).
    \end{multline*}
  \end{proof}
\end{prop}

\subsection{Homogeneous forms on $\con$}

We consider forms in $\Lambda^k(\con)$ whose coefficients are homogeneous of degree $r$, and show how these correspond to $k$-forms on the simplex $T=\{x\in\cO\mid s=1\}$. For the purposes of this section, the homogeneous forms are smooth but need not be polynomials, so $r$ need not be positive or an integer, though it will be both in practice.
\begin{dfn}
  Let $\Lambda^k_r(\con)$ denote the forms in $\Lambda^k(\con)$ whose coefficients are homogeneous of degree $r$. That is, $\alpha$ is in $\Lambda^k_r(\con)$ if it satisfies
  \begin{equation*}
    \alpha_{sx}(X_1,\dotsc,X_k):=s^r\alpha_{x}(X_1,\dotsc,X_k)
  \end{equation*}
  for all $x\in T$, $s>0$, and vectors $X_1,\dotsc,X_k$ in $\mb R^{n+1}$.
\end{dfn}

Note that $s$ has homogeneous degree one, $ds$ has homogeneous degree zero, and $X$ has homogeneous degree one in the sense that $X_{sx}=sX_x$. Hence, we observe that
\begin{equation*}
  \begin{tikzcd}[row sep=0]
    (ds\wedge)\colon\Lambda^k_r(\con)\arrow[r]&\Lambda^{k+1}_r(\con)\\
    i_X\colon\Lambda^k_r(\con)\arrow[r]&\Lambda^{k-1}_{r+1}(\con).
  \end{tikzcd}
\end{equation*}

We have results analogous to the ones in the previous subsection.
\begin{prop}\label{dsixexacthom}
  The following sequences are exact.
  \begin{equation*}
    \begin{tikzcd}[row sep=0]
      \Lambda_r^{k-1}(\con)\arrow[r,"ds\wedge"]&\Lambda_r^k(\con)\arrow[r,"ds\wedge"]&\Lambda_r^{k+1}(\con)\\
      \Lambda_{r+1}^{k-1}(\con)&\Lambda_r^k(\con)\arrow[l,"i_X" above]&\Lambda^{k+1}_{r-1}(\con)\arrow[l,"i_X" above]
    \end{tikzcd}
  \end{equation*}
  \begin{proof}
    The same proof as in Corollary \ref{dsixexact} applies, with the additional note that, if $\alpha\in\Lambda^k_r(\con)$, then $s^{-1}i_X\alpha\in\Lambda^{k-1}_r(\con)$, and $s^{-1}ds\wedge\alpha\in\Lambda^{k+1}_{r-1}(\con)$.
  \end{proof}
\end{prop}

Thus, we can make an analogous definition.
\begin{dfn}\label{vhdfn}
  Let
  \begin{align*}
    \Lambda^k_r(\con)^\perp&=ds\wedge\Lambda^{k-1}_r(\con)=\ker(ds\wedge)\cap\Lambda^k_r(\con),\\
    \Lambda^k_r(\con)^\top&=i_X\Lambda^{k+1}_{r-1}(\con)=\ker i_X\cap\Lambda^k_r(\con).
  \end{align*}
\end{dfn}

We then have analogous propositions.

\begin{prop}
  We have the following splitting of $\Lambda^k(\con)$.
  \begin{equation*}
    \Lambda^k_r(\con)=\Lambda^k_r(\con)^\perp\oplus\Lambda^k_r(\con)^\top.
  \end{equation*}
  This splitting is pointwise orthogonal with respect to the inner product $g$.
  \begin{proof}
    Using the splitting in Proposition \ref{dsixsplitting}, it remains to note that, if $\alpha\in\Lambda^k_r(\con)$, then its two components, $ds\wedge i_X(s^{-1}\alpha)$ and $i_X(ds\wedge s^{-1}\alpha)$, have homogeneous degree $r$, so they are in $\Lambda^k_r(\con)^\perp$ and $\Lambda^k_r(\con)^\top$, respectively.
  \end{proof}
\end{prop}

We make note of a basic fact.
\begin{prop}
  Multiplication by $s^{r'-r}$ is an isomorphism
  \begin{equation*}
    \Lambda^k_r(\con)\xrightarrow{s^{r'-r}}\Lambda^k_{r'}(\con).
  \end{equation*}
  that is compatible with the splitting into vertical and horizontal forms above.
  \begin{proof}
    Compatibility with the splitting follows from the fact that multiplication by $s^{r'-r}$ (or any scalar field) commutes with the tensorial operations $ds\wedge$ and $i_X$.
  \end{proof}
\end{prop}

\begin{prop}\label{vhisohom}
  The vertical homogeneous $(k+1)$-forms are isomorphic to the horizontal homogeneous $k$-forms via any of the following pairs of inverse maps.
  \begin{equation*}
    \begin{tikzcd}
      \Lambda^{k+1}_{r-1}(\con)^\perp \arrow[r, shift left, "i_X"] &\Lambda^k_r(\con)^\top \arrow[l, shift left, "s^{-1}ds\wedge"]&
      \Lambda^{k+1}_r(\con)^\perp \arrow[r, shift left, "s^{-1}i_X"] &\Lambda^k_r(\con)^\top \arrow[l, shift left, "ds\wedge"]\\
      \Lambda^{k+1}_{r-1}(\con)^\perp \arrow[r, shift left, "j_X"] &\Lambda^k_r(\con)^\top \arrow[l, shift left, "s^{-1}(\wedge ds)"]&
      \Lambda^{k+1}_r(\con)^\perp \arrow[r, shift left, "s^{-1}j_X"] &\Lambda^k_r(\con)^\top \arrow[l, shift left, "\wedge ds"]
    \end{tikzcd}
  \end{equation*}
  Here, $j_X$ denotes \emph{right} contraction with $X$.
  \begin{proof}
    Using Proposition \ref{vhiso}, we simplify need to verify that the above maps raise or lower homogeneous degree as stated.
  \end{proof}
\end{prop}

\subsection{Differential forms on $T$, revisited}

We now turn to forms on the simplex $T=\{x\in\cO\mid s=1\}$. Let $i\colon T\hookrightarrow\con$ denote the inclusion map, and hence $i^*\colon\Lambda^k(\con)\to\Lambda^k(T)$ denotes the restriction of forms.
\begin{prop}\label{homtprop}
  The restriction map from horizontal homogeneous $k$-forms on $\con$ to $k$-forms on the simplex $T$ is an isomorphism. That is,
  \begin{equation*}
    i^*\colon\Lambda^k_r(\con)^\top\xrightarrow\cong\Lambda^k(T).
  \end{equation*}
  \begin{proof}
    To show surjectivity, our task is, given a form $a\in\Lambda^k(T)$, to construct a homogeneous horizontal extension $\alpha\in\Lambda^k_r(\con)^\top$. For $x\in T$, $s>0$, and $X_1,\dotsc,X_k\in\mb R^{n+1}$, define
    \begin{equation*}
      \alpha_{sx}(X_1,\dotsc,X_k):=(s^r\pi_x^*a_{x})(X_1,\dotsc,X_k)
    \end{equation*}
    It is clear that $\alpha$ is homogeneous of degree $r$ and that $\alpha$ agrees with $a$ when $s=1$ and $X_1,\dotsc,X_k$ are tangent to $T$. To show that $i_X\alpha=0$, we use the fact that $X_{sx}=sX_x$ and the fact that $\pi_x^*a_x\in\Lambda^k_x(\con)^\top$, so $i_{X_x}(\pi_x^*a_x)=0$ for $x\in T$. For every $s>0$, we have
    \begin{equation*}
      \left(i_X\alpha\right)_{sx}=i_{X_{sx}}\alpha_{sx}=si_{X_x}(s^r\pi_x^*a_x)=s^{r+1}i_{X_x}(\pi_x^*a_x)=0.
    \end{equation*}
    Thus, $i_X\alpha=0$ on $\con$, so $\alpha\in\Lambda^k_r(\con)^\top$, as desired.

    To show injectivity, assume that $i^*\alpha=0$ for some $\alpha\in\Lambda^k_r(\con)^\top$. By Corollary \ref{htcor}, $i^*$ is injective on $\Lambda^k_x(\con)^\top$ for $x\in T$. Hence, $\alpha_x=0$ for any $x\in T$. But, since $\alpha$ is homogeneous, $\alpha_{sx}=s^r\alpha_x=0$ for any $s>0$. Thus $\alpha=0$ on all of $\con$.
  \end{proof}
\end{prop}
We let $h_r$ denote the extension map in the above proof. See Example \ref{extensioneg}.
\begin{dfn}
  Let $h_r\colon\Lambda^k(T)\to\Lambda^k_r(\con)^\top$ denote the map defined by
  \begin{equation*}
    (h_ra)_{sx}=s^r\pi_x^*a_x
  \end{equation*}
  for $x\in T$ and $s>0$.
\end{dfn}

With this definition, we can restate Proposition \ref{homtprop} as the following corollary.
\begin{cor}\label{homtcor}
  The space of forms on the simplex $T$ is isomorphic to the space of homogeneous horizontal forms on $\con$ via the horizontal homogeneous extension map $h_r$ and the restriction map $i^*$.
  \begin{equation*}
    \begin{tikzcd}
     \Lambda^k_r(\con)^\top \arrow[r, shift left, "i^*"] &  \Lambda^k(T) \arrow[l, shift left, "h_r"].
    \end{tikzcd}
  \end{equation*}
\end{cor}

Combining Corollary \ref{homtcor} with Proposition \ref{koszulcommutes} along with the observation that $X^\kappa$ has homogenenous degree one, we have the following corollary.

\begin{cor}\label{koszulcommutescor}
  The following diagram commutes
  \begin{equation*}
    \begin{tikzcd}
      \Lambda^{k+1}_{r-1}(\con)^\top \arrow[r, shift left, "i^*"] \arrow[d, "\kappa"]&\Lambda^{k+1}(T). \arrow[l, shift left, "h_{r-1}"] \arrow[d,"\kappa"] \\
      \Lambda^k_r(\con)^\top \arrow[r, shift left, "i^*"] &\Lambda^k(T). \arrow[l, shift left, "h_r"] 
    \end{tikzcd}
  \end{equation*}
\end{cor}

Combining Corollary \ref{homtcor} with Proposition \ref{vhisohom} and the fact that $i^*s^{-1}=1$, we have the following corollary. The motivation for the specific choice of isomorphism using $(\wedge ds)$ rather than $(ds\wedge)$ will become apparent in Theorem \ref{dcommutesvert}.

\begin{cor}\label{vertcor}
  The space of $k$-forms on the simplex $T$ is isomorphic to the space of homogeneous vertical $(k+1)$-forms on $\con$ via the following maps.
  \begin{equation*}
    \begin{tikzcd}[column sep=large]
      \Lambda^{k+1}_r(\con)^\perp \arrow[r, shift left, "i^*\circ j_X"] & \Lambda^k(T) \arrow[l, shift left, "(\wedge ds)\circ h_r"]
    \end{tikzcd}
  \end{equation*}
\end{cor}

In practice, to compute $h_r(a)$, it is easiest to construct an arbitrary homogeneous extension $\alpha'\in\Lambda^k_r(\con)$ of $a$ and then project $\alpha'$ to $\alpha\in\Lambda^k_r(\con)^\top$ via Proposition \ref{dsixs}.
\begin{prop}\label{homtohorhom}
  Let $\alpha'\in\Lambda^k_r(\con)$ and let $i^*\alpha'=a\in\Lambda^k(T)$. Then
  \begin{equation*}
    h_r(a)=s^{-1}i_X(ds\wedge\alpha')=\alpha'-s^{-1}ds\wedge i_X\alpha'.
  \end{equation*}
  \begin{proof}
    The equation
    \begin{equation*}
      s^{-1}i_X(ds\wedge\alpha')=\alpha'-s^{-1}ds\wedge i_X\alpha'
    \end{equation*}
    follows from Proposition \ref{dsixs}. Let $\alpha$ denote this expression, so we aim to show that $h_r(a)=\alpha$.

    First, note that $\alpha'$ having homogeneous degree $r$ implies that $\alpha$ has homogeneous degree $r$ as well. Next, observe that
    \begin{equation*}
      i_X\alpha=i_Xs^{-1}i_X(ds\wedge\alpha')=s^{-1}i_Xi_X(ds\wedge\alpha')=0.
    \end{equation*}
    Thus $\alpha\in\Lambda^k_r(\con)^\top$, so by Corollary \ref{homtcor}, $\alpha=h_r(i^*\alpha)$.

    Next, we see that the restriction of $\alpha$ is $a$.
    \begin{equation*}
      i^*\alpha=i^*\alpha'-i^*\left(s^{-1}ds\wedge i_X\alpha'\right)=a-i^*(ds)\wedge i^*\left(i_X\alpha'\right)=a
    \end{equation*}
    since $i^*(ds)=0$. Thus $\alpha=h_r(a)$, as desired.
  \end{proof}
\end{prop}

\begin{eg}\label{extensioneg}
  We provide some examples of homogeneous horizontal extensions for the two-dimensional simplex $s=x+y+z=1$ in the first octant of $\mb R^{3}$. We will use notation like $y\,dx$ both for forms on $T$ and for forms on $\cO$. We computed these extensions using Proposition \ref{homtohorhom}. One can also directly verify that the extensions below are homogeneous of the specified degree, are in the kernel of $i_X$, and have the specified restriction to $\Lambda^k(T)$.
  \begin{align*}
    h_4(x^4+3xy+y^3)&=x^4+3s^2xy+sy^3,\\
    h_2(x^4+3xy+y^3)&=s^{-2}x^4+3xy+s^{-1}y^3,\\
    h_1(dx)&=s\,dx-x\,ds,\\
    h_0(dx)&=dx-\tfrac xs\,ds,\\
    h_1(y\,dx)&=y\,dx-\tfrac{xy}s\,ds,\\
    h_1(y\,dx-x\,dy)&=y\,dx-x\,dy,\\
    h_1(dx\wedge dy)&=z\,dx\wedge dy+x\,dy\wedge dz+y\,dz\wedge dx
    .
  \end{align*}
  We see from these examples that, for scalar fields, the construction of the extension is effectively the same as the standard procedure for homogenization. But for higher degree forms, the situation is more interesting.

  One important thing to note is that, even though $y\,dx$ has polynomial coefficients of degree one, $h_1(y\,dx)$ does \emph{not} have polynomial coefficients. On the other hand, $h_1(y\,dx-x\,dy)$ does have polynomial coefficients. As we will see in Theorem \ref{prmthm}, this behavior occurs because $y\,dx$ is not in $\mc P_1^-\Lambda^1(T)$, whereas $y\,dx-x\,dy$ is.
\end{eg}

\subsection{Vanishing tangential trace}
The notation $\mr{\Lambda}^k(T)$ denotes the forms on $T$ whose restriction as differential forms to $\partial T$ is zero. We define the corresponding space of forms on $\con$.

\begin{dfn}
  Let $\partial\con$ denote the nonzero points of $\partial\cO$, and let $k\colon\partial\con\to\con$ be the inclusion. Let $\mr{\Lambda}^k(\con)$ denote the forms whose restriction to the $\partial\con$ is zero, that is, those forms $\alpha\in\Lambda^k(\con)$ such that $k^*\alpha=0$.
\end{dfn}

\begin{prop}\label{hornotrace}
  The isomorphism in Corollary \ref{homtcor} preserves the property of vanishing trace. That is, it restricts to the isomorphism
  \begin{equation*}
    \begin{tikzcd}
     \mr\Lambda^k_r(\con)^\top \arrow[r, shift left, "i^*"] &  \mr\Lambda^k(T) \arrow[l, shift left, "h_r"].
    \end{tikzcd}
  \end{equation*}
  \begin{proof}
    The forward direction is clear because $\partial T\subset\partial\con$. Indeed, if $\alpha\in\Lambda^k_r(\con)$ vanishes when restricted to $\partial\con$, then it vanishes when further restricted to $\partial T$.

    For the reverse direction, let $\alpha\in\Lambda^k_r(\con)^\top$, and assume that $i^*\alpha$ has vanishing trace on $\partial T$. We aim to show that $\alpha$ has vanishing trace on $\partial\con$.

    Let $x\in\partial T$, and let $X_1,\dotsc,X_k$ be arbitrary tangent vectors in $T_x\partial\con$. Crucially, $X$ is also tangent to $\partial\con$, and $T_x\partial\con$ is spanned by $T_x\partial T$ and $X$. Hence, we can write $X_i=Y_i+c_iX$, where $Y_i\in T_x\partial T$. The assumption that $i^*\alpha$ has vanishing trace on $\partial T$ tells us that $\alpha_x(Y_1,\dotsc,Y_k)=0$. As in the proof of Proposition \ref{htiso}, this fact, along with the fact that $i_X\alpha=0$, tell us that $\alpha_x(X_1,\dotsc,X_k)=0$. In other words, the restriction of $\alpha$ to $\partial\con$ vanishes at the point $x$. But any point in $\partial\con$ is of the form $sx$ where $s>0$ and $x\in\partial T$, so, by homogeneity of $\alpha$, the restriction of $\alpha$ to $\partial\con$ vanishes at all points of $\con$.
  \end{proof}
\end{prop}

\begin{prop}\label{vertnotrace}
  The isomorphism in Corollary \ref{vertcor} preserves the property of vanishing trace. That is, it restricts to the isomorphism
  \begin{equation*}
    \begin{tikzcd}[column sep=large]
      \mr\Lambda^{k+1}_r(\con)^\perp \arrow[r, shift left, "i^*\circ j_X"] & \mr\Lambda^k(T). \arrow[l, shift left, "(\wedge ds)\circ h_r"]
    \end{tikzcd}
  \end{equation*}
  \begin{proof}
    Given Proposition \ref{hornotrace}, it suffices to show that the bottom right isomorphism from Proposition \ref{vhisohom} preserves the property of vanishing trace. That is, it suffices to show the isomorphism from Proposition \ref{vhisohom} restricts to
    \begin{equation*}
      \begin{tikzcd}
        \mr\Lambda^{k+1}_r(\con)^\perp \arrow[r, shift left, "s^{-1}j_X"] &\mr\Lambda^k_r(\con)^\top. \arrow[l, shift left, "\wedge ds"]
      \end{tikzcd}
    \end{equation*}
    Pullbacks commute with scalar multiplication and wedge multiplication, and, because $X$ is tangent to $\partial\con$, the pullback $k^*$ commutes with $j_X$. Thus, if $k^*\alpha=0$, then $k^*(s^{-1}j_X\alpha)=s^{-1}j_Xk^*\alpha=0$, and if $k^*\beta=0$, then $k^*(\beta\wedge ds)=k^*\beta\wedge k^*(ds)=0$.
  \end{proof}
\end{prop}

\section{The spaces $\mbf P_r\Lambda^k$ and $\mbf P_r^-\Lambda^k$}\label{prprmsec}
When we restrict to forms on $\cO$ whose coefficients are polynomials, the spaces of vertical and horizontal homogeneous forms naturally give the $\mc P_r$ and $\mc P_r^-$ subspaces of Arnold, Falk, and Winther \cite{afw06,afw10} via the isomorphisms in Corollaries \ref{homtcor} and \ref{vertcor}.

\begin{dfn}
  Define the spaces of homogeneous polynomial differential forms on $\cO$ as follows.
  \begin{align*}
    \mbf H_r\Lambda^k(\cO)&=\{\text{polynomial differential forms}\}\cap\Lambda^k_r(\con),\\
    \mbf P_r\Lambda^k(\cO)&=\{\text{polynomial differential forms}\}\cap\Lambda^{k+1}_r(\con)^\perp,\\
    \mbf P_r^-\Lambda^k(\cO)&=\{\text{polynomial differential forms}\}\cap\Lambda^k_r(\con)^\top.
  \end{align*}
\end{dfn}
In other words, the $\mbf H$ spaces contain all homogeneous polynomial differential forms, then $\mbf P$ spaces contain the vertical homogeneous polynomial differential forms, and the $\mbf P^-$ spaces contain the horizontal homogeneous polynomial differential forms. However, note that, in light of Corollary \ref{vertcor}, the space $\mbf P_r\Lambda^k(\cO)$ is a space of $(k+1)$-forms.

Note also that, although $\Lambda^k_r(\con)=\Lambda^k_r(\con)^\perp\oplus\Lambda^k_r(\con)^\top$, when we restrict to polynomials, $\mbf H_r\Lambda^k(\cO)$ is in general larger than $\mbf P_r\Lambda^{k-1}(\cO)\oplus\mbf P_r^-\Lambda^k(\cO)$. In other words, there are homogeneous polynomial differential forms in $\Lambda^k_r(\con)$ whose projections to the spaces $\Lambda^k_r(\con)^\perp$ and $\Lambda^k_r(\con)^\top$ are not polynomials.

Finally, although $\Lambda^k_r(\con)$ excludes the origin, any polynomial defined on $\con$ is also defined on $\cO$, so omitting the $\neq0$ is acceptable.

\begin{eg}
    We can compute these spaces either from the definition or by making use of the exact sequences given by Proposition \ref{polyexact}. With $n=2$, we have
  \begin{align*}
    \mbf P_1^-\Lambda^0(\cO)&=\Span\{x,y,z\},&
    \mbf P_1\Lambda^0(\cO)&=\Span\{x\,ds, y\,ds, z\,ds\},\\
    \begin{split}
    \mbf P_1^-\Lambda^1(\cO)&=\Span\{y\,dx-x\,dy,\\
    &\phantom{{}=\Span\{} z\,dy-y\,dz,\\
    &\phantom{{}=\Span\{} x\,dz-z\,dx\},
    \end{split}&
    \begin{split}
      \mbf P_1\Lambda^1(\cO)&=\Span\{y\,dx\wedge ds,x\,dy\wedge ds,\\
      &\phantom{{}=\Span\{} z\,dy\wedge ds,y\,dz\wedge ds,\\
      &\phantom{{}=\Span\{} x\,dz\wedge ds,z\,dx\wedge ds\},
    \end{split}\\
    \begin{split}
      \mbf P_1^-\Lambda^2(\cO)&=\Span\{x\,dy\wedge dz\\
      &\qquad\quad {}+y\,dz\wedge dx\\
      &\qquad\quad {}+z\,dx\wedge dy\},
    \end{split}&
    \begin{split}
      \mbf P_1\Lambda^2(\cO)&=\Span\{x\,dx\wedge dy\wedge dz,\\
      &\phantom{{}=\Span\{} y\,dx\wedge dy\wedge dz,\\
      &\phantom{{}=\Span\{} z\,dx\wedge dy\wedge dz\}.
    \end{split}
  \end{align*}
\end{eg}


By definition, we have
\begin{align*}
  \mbf P_r\Lambda^k(\cO)&=\mbf H_r\Lambda^{k+1}(\cO)\cap\ker(ds\wedge),\\
  \mbf P_r^-\Lambda^k(\cO)&=\mbf H_r\Lambda^k(\cO)\cap\ker i_X.
\end{align*}
We would like to express the $\mbf P_r$ and $\mbf P_r^-$ spaces as the image of $\mbf H_r$ under the maps $ds\wedge$ and $i_X$, similarly to Definition \ref{vhdfn}. However, we cannot simply use Proposition \ref{dsixexacthom} because multiplication by $s^{-1}$ used in the proof does not generally result in a polynomial. Nonetheless, an analogous result holds.

\begin{prop}\label{polyexact}
  The following sequences are exact, except in the case $r=k=0$, in which case only the first sequence is exact.
  \begin{equation*}
    \begin{tikzcd}[row sep=0]
      \mbf H_r\Lambda^{k-1}(\cO)\arrow[r,"ds\wedge"]&\mbf H_r\Lambda^k(\cO)\arrow[r,"ds\wedge"]&\mbf H_r\Lambda^{k+1}(\cO)\\
      \mbf H_{r+1}\Lambda^{k-1}(\cO)&\mbf H_r\Lambda^k(\cO)\arrow[l,"i_X" above]&\mbf H_{r-1}\Lambda^{k+1}(\cO)\arrow[l,"i_X" above]
    \end{tikzcd}
  \end{equation*}
  In other words, after adjusting the index $k$,
  \begin{align*}
    \mbf P_r\Lambda^k(\cO)&=ds\wedge \mbf H_r\Lambda^k(\cO),\\
    \mbf P_r^-\Lambda^k(\cO)&=i_X\mbf H_{r-1}\Lambda^{k+1}(\cO).
  \end{align*}

  \begin{proof}
    Given Proposition \ref{dsixexacthom}, we know that if $\alpha\in\mbf H_r\Lambda^k(\cO)$ satisfies $ds\wedge\alpha=0$, then it is in the image of $ds\wedge$, but what we need to show is that $\alpha=ds\wedge\beta$ where $\beta$ is a \emph{polynomial} in $\mbf H_r\Lambda^{k-1}(\cO)$, and similarly for the $i_X$ operator.

    By Lemma \ref{dslemma}, if $ds\wedge\alpha=0$, then $\alpha=ds\wedge i_{\nabla s}\left(\frac1{n+1}\alpha\right)$. Since $\nabla s$ is a constant vector field, $i_{\nabla s}\left(\frac1{n+1}\alpha\right)\in\mbf H_r\Lambda^{k-1}(\cO)$, as desired.

    Similarly, by Lemma \ref{ixlemma}, if $i_X\alpha=0$, then $\alpha=i_X\left(\frac1{r+k}d\alpha\right)$. The exterior derivative decreases polynomial degree and increases form degree, so $\frac1{r+k}d\alpha\in\mbf H_{r-1}\Lambda^{k+1}(\cO)$, as desired.
  \end{proof}
\end{prop}

Note that the latter argument fails when $r=k=0$, and, indeed $\mbf H_0\Lambda^0(\cO)$ is a one-dimensional space consisting of constant scalar fields, whereas $\mbf H_1\Lambda^{-1}(\cO)$ and $\mbf H_{-1}\Lambda^1(\cO)$ are zero.

\begin{lemma}\label{dslemma}
  We have
  \begin{equation*}
    (ds\wedge)\circ i_{\nabla s}+i_{\nabla s}\circ(ds\wedge)=n+1.
  \end{equation*}
  That is, for any $\alpha\in\Lambda^k(\mbf O)$,
  \begin{equation*}
    ds\wedge i_{\nabla s}\alpha+i_{\nabla s}(ds\wedge\alpha)=(n+1)\alpha.
  \end{equation*}
  \begin{proof}
    The proof is identical to that of Proposition \ref{dsixs}, except that we have $\nabla s=\sum_{i=1}^{n+1}\pp{x_i}$, so
    \begin{equation*}
      i_{\nabla s}(ds)=n+1.\qedhere
    \end{equation*}
  \end{proof}
\end{lemma}

\begin{lemma}\label{ixlemma}
  We have
  \begin{equation*}
    d\circ i_X+i_X\circ d = \mf L_X = r+k,
  \end{equation*}
  where $\mf L_X$ denotes the Lie derivative. That is, for $\alpha\in\mbf H_r\Lambda^k(\cO)$,
  \begin{equation*}
    d(i_X\alpha)+i_X(d\alpha)=\mf L_X\alpha=(r+k)\alpha.
  \end{equation*}
  \begin{proof}
    The statement that $d\circ i_X+i_X\circ d=\mf L_X$ is Cartan's formula. For the statement that for $\alpha\in\mbf H_r\Lambda^k(\cO)$, we have $\mf L_X\alpha=(r+k)\alpha$, see a proof of an analogous statement in \cite[Theorem 3.1]{afw06}. Alternatively, observe that
    \begin{equation*}
      \mf L_Xx_i=d(i_Xx_i)+i_X(dx_i)=0+x_i=x_i.
    \end{equation*}
    Algebraic properties of the Lie derivative tell us that the Lie derivative commutes with the exterior differential and satisfies the Leibniz rule with respect to the wedge product. That is, $\mf L_X(d\beta)=d(\mf L_X\beta)$ and $\mf L_X(\beta\wedge\gamma)=\mf L_X\beta\wedge\gamma+\beta\wedge\mf L_X\gamma$, notably including the case where $\beta$ is a scalar field, that is, a $0$-form. These two properties, along with $\mf L_Xx_i=x_i$, suffice to show that
    \begin{multline*}
      \mf L_X\left(x_{i_1}x_{i_2}\dotsm x_{i_r}dx_{j_1}\wedge\dotsb\wedge dx_{j_k}\right)\\=(r+k)\left(x_{i_1}x_{i_2}\dotsm x_{i_r}dx_{j_1}\wedge\dotsb\wedge dx_{j_k}\right).
    \end{multline*}
    The claim follows by linearity.
  \end{proof}
\end{lemma}

We recall the definitions of the $\mc P_r$ and $\mc P_r^-$ spaces of Arnold, Falk, and Winther \cite{afw06,afw10}.

\begin{dfn}
  The notation $\mc P_r\Lambda^k(T)$ denotes the space $k$-forms on the $n$-simplex $T$ whose coefficients are polynomials of degree at most $r$. The notation $\mc P_r^-\Lambda^k(T)$ denotes the subset of $\mc P_r\Lambda^k(T)$ defined by
  \begin{equation*}
    \mc P_r^-\Lambda^k(T):=\mc P_{r-1}\Lambda^k(T)+\kappa\mc P_{r-1}\Lambda^{k+1}(T),
  \end{equation*}
  where $\kappa$ denotes the Koszul operator. (See Definitions \ref{koszulvfdfn} and \ref{koszuldfn}.)
\end{dfn}

We now prove that the $\mbf P_r$ and $\mbf P_r^-$ spaces are isomorphic to the $\mc P_r$ and $\mc P_r^-$ spaces via the isomorphisms in Corollaries \ref{homtcor} and \ref{vertcor}. The main subtlety in the proof is that the horizontal homogeneous extension $h_r(a)$ of $a\in\Lambda^k(T)$ need not have polynomial coefficients even if $a$ has coefficients that are polynomials of degree $r$, due to the fact that the horizontal subspace of forms varies with $x\in T$. However, if we remove the horizontal condition, then it is easy to construct a homogeneous polynomial extension of $a$.

\begin{lemma}
  Any differential form on $T$ with polynomial coefficients of degree $r$ can be extended to a differential form on $\cO$ with homogeneous \emph{polynomial} coefficients of degree $r$. That is, given $a\in\mc P_r\Lambda^k(T)$, there exists an $\alpha'\in\mbf H_r\Lambda^k(\cO)$ such that $i^*\alpha'=a$.
  \begin{proof}
    One way to proceed is to replicate our construction of $h_r$ to construct an extension $h_r'(a)$ using the standard orthogonal projection to the tangent space instead of $\pi_x$ and using $\nabla s$ instead of $X$. However, none of our constructions need a specific choice of homogeneous polynomial extension. Thus, it suffices to construct an arbitrary polynomial extension $\alpha'$; there is no need to impose the additional constraint that $i_{\nabla s}\alpha'=0$.

    Hence, we instead proceed by simply showing that $i^*\mbf H_r\Lambda^k(\cO)$ spans $\mc P_r\Lambda^k(T)$. We start with zero-forms, where we essentially do the standard procedure of homogenizing/dehomogenizing polynomials, except with $s=1$ instead of $x_{n+1}=1$. First, consider $i^*\mbf H_1\Lambda^0(\cO)$. Observe that $i^*x_i=x_i$ for $1\le i\le n+1$, and note that $i^*(x_1+\dotsb+x_{n+1})=i^*s=1$ on $T$. Thus, these $n+1$ functions span all linear functions on $T$, so $i^*\mbf H_1\Lambda^0(\cO)=\mc P_1\Lambda^0(T)$. Next, consider a product of linear factors $p=l_1\dotsm l_r$, where the $l_j$ are in $\mbf H_1\Lambda^0(\cO)$. Then $p\in\mbf H_r\Lambda^0(\cO)$, and $i^*p=(i^*l_1)\dotsm(i^*l_r)$ is a product of $r$ functions in $\mc P_1\Lambda^0(T)$. Products of $r$ functions in $\mc P_1\Lambda^0(T)$ span $\mc P_r\Lambda^0(T)$, so $i^*\mbf H_r\Lambda^0(\cO)=\mc P_r\Lambda^0(\cO)$. (In fact, the restriction map $i^*\colon\mbf H_r\Lambda^0(\cO)\to\mc P_r\Lambda^0(T)$ is an isomorphism on $0$-forms.)

    Next, consider $i^*\mbf H_0\Lambda^1(\cO)$. Observe that the $n+1$ one-forms $i^*dx_i$ for $1\le i\le n+1$ span $\mc P_0\Lambda^1(T)$ with the single relation $i^*dx_1+\dotsb+i^*dx_{n+1}=i^*(ds)=0$. Thus, $i^*\mbf H_0\Lambda^1(\cO)=\mc P_0\Lambda^1(T)$. Now consider $\alpha=p\,\theta_1\wedge\dotsb\wedge\theta_k$, where $p$ is a scalar polynomial in $\mbf H_r\Lambda^0(\cO)$ and the $\theta_j$ are constant one-forms in $\mbf H_0\Lambda^1(\cO)$. Then $\alpha\in\mbf H_r\Lambda^k(\cO)$, and $i^*\alpha=i^*p\,(i^*\theta_1)\wedge\dotsb\wedge(i^*\theta_k)$. The wedge products of $k$ one-forms in $\mc P_0\Lambda^1(T)$ span $\mc P_0\Lambda^k(T)$, and products of polynomials in $\mc P_r\Lambda^0(T)$ and constant $k$-forms in $\mc P_0\Lambda^k(T)$ span $\mc P_r\Lambda^k(T)$. Hence, $i^*\mbf H_r\Lambda^k(\cO)=\mc P_r\Lambda^k(T)$, as desired.
  \end{proof}
\end{lemma}

\begin{eg}
  For example, $dx+x\,dy\in\mc P_1\Lambda^1(T)$ can be extended to the homogeneous degree one polynomial $s\,dx+x\,dy\in\mbf H_1\Lambda^1(\cO)$. Since, when restricted to $T$, $dx+dy+dz=ds=0$, we could, alternatively, have extended $dx+x\,dy$ to $-s\,dy-s\,dz+x\,dy\in\mbf H_1\Lambda^1(\cO)$.
\end{eg}

\begin{thm}\label{prmthm}
  Let $r\ge1$. The isomorphism of Corollary \ref{homtcor} restricts to the following isomorphism between polynomial spaces
  \begin{equation*}
    \begin{tikzcd}
      \mbf P_r^-\Lambda^k(\cO) \arrow[r, shift left, "i^*"] & \mc P_r^-\Lambda^k(T) \arrow[l, shift left, "h_r"]
    \end{tikzcd}
  \end{equation*}
  \begin{proof}
    Recall that
    \begin{equation*}
      \mc P_r^-\Lambda^k(T)=\mc P_{r-1}\Lambda^k(T)+\kappa\mc P_{r-1}\Lambda^{k+1}(T).
    \end{equation*}
    
    Let $\alpha\in\mbf P_r^-\Lambda^k(\cO)$. We aim to show that $i^*\alpha\in\mc P_r^-\Lambda^k(T)$. By Proposition \ref{polyexact}, $\alpha=i_X\beta$ for some $\beta\in\mbf H_{r-1}\Lambda^{k+1}(\cO)$. Using Corollary \ref{koszulcommutescor}, Definitions \ref{koszulvfdfn} and \ref{koszuldfn}, and the fact that $i^*s=1$, we find that
    \begin{equation*}
      i^*\alpha=i^*i_X\beta=i^*\left(\tfrac1{n+1}si_{\nabla s}\beta\right)+i^*(\kappa\beta)=\tfrac1{n+1}i^*\left(i_{\nabla s}\beta\right)+\kappa i^*\beta.
    \end{equation*}
    Because $\nabla s$ is a constant vector field, $i_{\nabla s}\beta$ is a polynomial in $\mbf H_{r-1}\Lambda^k(\cO)$, and hence its restriction to $T$ is a polynomial in $\mc P_{r-1}\Lambda^k(T)$. Likewise, $i^*\beta$ is a polynomial in $\mc P_{r-1}\Lambda^{k+1}(T)$. Hence
    \begin{equation*}
      i^*\alpha\in\mc P_{r-1}\Lambda^k(T)+\kappa\mc P_{r-1}\Lambda^{k+1}(T),
    \end{equation*}
    as desired.

    Conversely, let $a\in\mc P_r^-\Lambda^k(T)$. We must show that $h_r(a)$ is a polynomial. If $a\in\mc P_{r-1}\Lambda^k$, then we can let $\beta\in\mbf H_{r-1}\Lambda^k(\cO)$ be an arbitrary homogeneous polynomial extension of $a$ of degree $r-1$. Then $\alpha'=s\beta\in\mbf H_r\Lambda^k(\cO)$ is a homogeneous polynomial extension of $a$ of degree $r$. Consequently, by Proposition \ref{homtohorhom},
    \begin{equation*}
      h_r(a)=s^{-1}i_X(ds\wedge\alpha')=s^{-1}i_X(ds\wedge s\beta)=i_X(ds\wedge\beta),
    \end{equation*}
    which is a polynomial.

    Meanwhile, if $a=\kappa b$ for $b\in\mc P_{r-1}\Lambda^{k+1}$, then let $\beta\in\mbf H_{r-1}\Lambda^{k+1}(\cO)$ be an arbitrary homogeneous polynomial extension of $b$. Then, by Corollary \ref{koszulcommutescor}, Definitions \ref{koszulvfdfn} and \ref{koszuldfn}, and Proposition \ref{homtohorhom}, we find that
    \begin{multline*}
      h_r(a)=h_r(\kappa b)=\kappa\left(h_{r-1}(b)\right)=\left(i_X-\tfrac1{n+1}si_{\nabla s}\right)\left(s^{-1}i_X(ds\wedge\beta)\right)\\
      =-\tfrac1{n+1}i_{\nabla s}i_X(ds\wedge\beta)=\tfrac1{n+1}i_Xi_{\nabla s}(ds\wedge\beta),
    \end{multline*}
    since $i_Xi_X=0$ and multiplication by the scalar field $s$ commutes with contraction. Thus, $h_r(a)$ is indeed a polynomial.
  \end{proof}
\end{thm}

In fact, the above argument works perfectly well even when $r=0$, as long as $k\neq0$. However, in that case, both sides of the isomorphism are zero. We have a similar claim for the vertical differential forms.

\begin{thm}\label{prthm}
  The isomorphism in Corollary \ref{vertcor} restricts to the following isomorphism between polynomial spaces.
  \begin{equation*}
    \begin{tikzcd}[column sep=large]
      \mbf P_r\Lambda^k(\cO) \arrow[r, shift left, "i^*\circ j_X"] & \mc P_r\Lambda^k(T) \arrow[l, shift left, "(\wedge ds)\circ h_r"]
    \end{tikzcd}
  \end{equation*}
  \begin{proof}
    Let $\alpha\in\mbf P_r\Lambda^k(\cO)$. We aim to show that $i^*\left(j_X\alpha\right)$ is a polynomial of degree $r$ in $\mc P_r(T)$. Using Proposition \ref{polyexact}, we see that $\alpha=\beta\wedge ds$ for some $\beta\in\mbf H_r\Lambda^k(\cO)$. Then, by Proposition \ref{dsixs},
    \begin{equation*}
      i^*\left(j_X\alpha\right)=i^*\left(j_X(\beta\wedge ds)\right)=i^*\left(s\beta-j_X\beta\wedge ds\right)=i^*\beta,
    \end{equation*}
    since $i^*s=1$ and $i^*(ds)=0$. Since $\beta\in\mbf H_r\Lambda^k(\cO)$, we know that $i^*\beta\in\mc P_r\Lambda^k(T)$, so, indeed, $i^*\left(j_X\alpha\right)$ is a polynomial differential form of degree $r$.

    Conversely, let $a\in\mc P_r\Lambda^k(T)$. It is clear that $h_r(a)\wedge ds\in\Lambda^{k+1}_r(\con)^\perp$, so we must show that $h_r(a)\wedge ds$ is a polynomial. Let $\alpha'\in\mbf H_r\Lambda^k(\cO)$ be an arbitrary polynomial extension of $a$. By Proposition \ref{homtohorhom},
    \begin{equation*}
      h_r(a)\wedge ds=\left(\alpha'-s^{-1}ds\wedge i_X\alpha'\right)\wedge ds=\alpha'\wedge ds,
    \end{equation*}
    because $ds\wedge\omega\wedge ds=0$ for any form $\omega$. Since $\alpha'$ is a polynomial, so is $\alpha'\wedge ds$. Hence, $h_r(a)\wedge ds\in\mbf H_r\Lambda^{k+1}(\cO)$, as desired.
  \end{proof}
\end{thm}

Note that, in the proof, the specific choice of extension $h_r(a)$ was not relevant; any other extension gives the same isomorphism map. We could have, for instance, have used the standard metric instead of $g$, giving us $j_{\frac1{n+1}\nabla s}$ instead of $j_X$. In contrast, for Theorem \ref{prmthm}, the use of the nonstandard metric $g$ and the resulting choice of extension $h_r$ is essential.

\subsection{Basic properties}
The $\mc P_r$ and $\mc P_r^-$ spaces have additional structure, in particular, exterior differentiation and pairing via integration. The remainder of the article will describe these structures in terms of the $\mbf P_r$ and $\mbf P_r^-$ spaces. We begin, however, with two basic strucutres: the subspaces $\mr{\mc P}_r$ and $\mr{\mc P}_r^-$ of forms with vanishing tangential boundary trace, and the inclusions $\mc P_{r-1}\hookrightarrow\mc P_r^-\hookrightarrow\mc P_r$.

\begin{dfn}
  Let $\mr{\mbf H}_r\Lambda^k(\cO)$ denote those $k$-forms whose restriction to the boundary of $\cO$ vanishes. That is, if $k\colon\partial\cO\hookrightarrow\cO$ is the inclusion, then
  \begin{equation*}
    \mr{\mbf H}_r\Lambda^k(\cO)=\{\alpha\in\mbf H_r\Lambda^k(\cO)\mid k^*\alpha=0\}.
  \end{equation*}
  Define the spaces $\mr{\mbf P}_r\Lambda^k(\cO)$ and $\mr{\mbf P}_r^-\Lambda^k(\cO)$ likewise.
\end{dfn}

\begin{prop}\label{prnotrace}
  The isomorphisms in Theorems \ref{prmthm} and \ref{prthm} preserve the property of vanishing trace. That is, we have the isomorphisms
  \begin{equation*}
    \begin{tikzcd}[column sep=large,row sep=small]
      \mr{\mbf P}_r^-\Lambda^k(\cO) \arrow[r, shift left, "i^*"] & \mr{\mc P}_r^-\Lambda^k(T) \arrow[l, shift left, "h_r"]\\
      \mr{\mbf P}_r\Lambda^k(\cO) \arrow[r, shift left, "i^*\circ j_X"] & \mr{\mc P}_r\Lambda^k(T) \arrow[l, shift left, "(\wedge ds)\circ h_r"]
    \end{tikzcd}
  \end{equation*}
  \begin{proof}
    We combine Theorems \ref{prmthm} and \ref{prthm} with Propositions \ref{hornotrace} and \ref{vertnotrace}.
  \end{proof}
\end{prop}

We move on to the inclusions $\mc P_{r-1}\hookrightarrow\mc P_r^-\hookrightarrow\mc P_r$.
\begin{prop}
  The maps $j_X$ and $\wedge ds$ correspond to the inclusion maps $\mc P_{r-1}\Lambda^k(T)\hookrightarrow\mc P_r^-\Lambda^k(T)\hookrightarrow\mc P_r\Lambda^k(T)$ as expressed in the following commutative diagram.
  \begin{equation*}
    \begin{tikzcd}[column sep=large]
       \mbf P_{r-1}\Lambda^k(\cO) \arrow[d, "j_X"] \arrow[r, shift left, "i^*\circ j_X"] & \mc P_{r-1}\Lambda^k(T) \arrow[d, hook] \arrow[l, shift left, "(\wedge ds)\circ h_{r-1}"]\\
       \mbf P_r^-\Lambda^k(\cO) \arrow[d,"\wedge ds"]\arrow[r, shift left, "i^*"]&\mc P_r^-\Lambda^k(T) \arrow[d,hook] \arrow[l, shift left, "h_r"]\\
       \mbf P_r\Lambda^k(\cO) \arrow[r, shift left, "i^*\circ j_X"]&\mc P_r\Lambda^k(T) \arrow[l, shift left, "(\wedge ds)\circ h_r"]
    \end{tikzcd}
  \end{equation*}
  \begin{proof}
    The top square in the diagram commutes by definition; it amounts to verifying that $i^*(j_X\alpha)=(i^*\circ j_X)\alpha$.

    For the bottom square, let $\alpha\in\mbf P_r^-\Lambda^k(\cO)$. We are tasked with verifying that $i^*(j_X(\alpha\wedge ds))=i^*\alpha$. Proposition \ref{dsixs}, along with the fact that $j_X\alpha=(-1)^{k-1}i_X\alpha=0$, tells us that
    \begin{equation*}
      i^*(j_X(\alpha\wedge ds))=i^*(s\alpha)=i^*\alpha,
    \end{equation*}
    as desired.
  \end{proof}
\end{prop}
A notable consequence of the above analysis is that the induced map $\mbf P_{r-1}\Lambda^k(\cO)\to\mbf P_r\Lambda^k(\cO)$ is simply multiplication by $s$, and likewise so is the induced map $\mbf P_{r-1}^-\Lambda^k(\cO)\to\mbf P_r^-\Lambda^k(\cO)$.

\section{Differentiation}\label{diffsec}
A key feature of the $\mc P_r$ and $\mc P_r^-$ spaces is their behavior with respect to the exterior derivative $d$. In this section, we develop the corresponding operation for the $\mbf P_r$ and $\mbf P_r^-$ spaces.

\begin{dfn}
  Define $\mbf d\colon\Lambda^k(\con)\to\Lambda^{k+1}(\con)$ by
  \begin{equation*}
    \mbf d = d - s^{-1}ds\wedge\mf L_X.
  \end{equation*}
  That is, for $\alpha\in\Lambda^k(\con)$,
  \begin{equation*}
    \mbf d\alpha=d\alpha -s^{-1}ds\wedge\mf L_X\alpha,
  \end{equation*}
  where $\mf L_X$ denotes the Lie derivative.
\end{dfn}

In practice, we will work with differential forms with homogeneous coefficients $\Lambda^k_r(\con)$, and in this setting the operator $\mbf d$ can be expressed in a simpler manner.
\begin{prop}\label{dhomog}
  Restricted to homogeneous forms, the operator $\mbf d\colon\Lambda^k_r(\con)\to\Lambda^{k+1}_{r-1}(\con)$ has the formula
  \begin{equation*}
    \mbf d = s^{r+k}\circ d\circ s^{-(r+k)}.
  \end{equation*}
  That is, if $\alpha\in\Lambda^k_r(\con)$, then
  \begin{equation*}
    \mbf d\alpha=s^{r+k}d\left(s^{-(r+k)}\alpha\right).
  \end{equation*}
  \begin{proof}
    By Lemma \ref{ixlemma}, $\mf L_X\alpha=(r+k)\alpha$. We thus compute that
    \begin{multline*}
      s^{r+k}d\left(s^{-(r+k)}\alpha\right)=s^{r+k}\left(s^{-(r+k)}d\alpha-(r+k)s^{-(r+k)-1}ds\wedge\alpha\right)\\
      =d\alpha-(r+k)s^{-1}ds\wedge\alpha=d\alpha-s^{-1}ds\wedge\mf L_X\alpha=\mbf d\alpha,
    \end{multline*}
    as desired.
  \end{proof}
\end{prop}
We can think of $r+k$ as the \emph{total degree} of $\alpha$, that is, the sum of the homogeneous degree of the coefficients and the degree of the form. One way to interpret the above result is that if $\alpha=\beta\,s^{r+k}$, where $\beta$ has total degree zero, then $\mbf d\alpha=d\beta\,s^{r+k}$. In other words, to apply $\mbf d$ to $\alpha$, we multiply $\alpha$ by a power of $s$ to obtain a corresponding form $\beta$ of total degree zero, then apply $d$, and then multiply by a power of $s$ to once again have a form of total degree $r+k$.

The operator $\mbf d$ satisfies many of the same properties as $d$.
\begin{prop}\label{ddzero}
  We have $\mbf d\circ\mbf d=0$.
  \begin{proof}[Proof for homogeneous forms]
    On homogeneous forms, if $\alpha\in\Lambda^k_r(\con)$, then $d\alpha\in\Lambda^{k+1}_{r-1}(\con)$. Notably, $(r-1)+(k+1)=r+k$, and so
    \begin{multline*}
      (\mbf d\circ \mbf d)(\alpha)=\left(s^{r-1+k+1}\circ d\circ s^{-(r-1+k+1)}\circ s^{r+k}\circ d\circ s^{-(r+k)}\right)(\alpha)\\
      =\left(s^{r+k}\circ d\circ d\circ s^{-(r+k)}\right)(\alpha)=0.
    \end{multline*}

    Alternatively, using the preceding remark, we write $\alpha=\beta\,s^{r+k}$, so $\beta$ has total degree zero. Then $\mbf d\alpha=d\beta\,s^{r+k}$. Since $d\beta$ also has total degree zero, $\mbf d(\mbf d\alpha)=d(d\beta)\,s^{r+k}=0$.
  \end{proof}
  \begin{proof}[Proof for general forms]
    In practice, we will only need the statement for homogeneous forms, but, for completeness, general forms $\alpha\in\Lambda^k(\con)$ can be handled in one of two ways. The first way is to decompose $\alpha$ into homogeneous components.
    \begin{equation*}
      \alpha=\int_{-\infty}^\infty \alpha_{r+k}s^{r+k}\,dr
    \end{equation*}
    where $\alpha_{r+k}$ has total degree zero. Thus, as discused above, $\mbf d(\alpha_{r+k}s^{r+k})=d\alpha_{r+k}s^{r+k}$, and so we have
    \begin{equation*}
      \mbf d\alpha=\int_{-\infty}^\infty d\alpha_{r+k}s^{r+k}\,dr.
    \end{equation*}
    Because $d\alpha_{r+k}$ has total degree zero, we have
    \begin{equation*}
      \mbf d(\mbf d\alpha)=\int_{-\infty}^\infty d(d\alpha_{r+k})s^{r+k}\,dr=0
    \end{equation*}
    because $d(d\alpha_{r+k})=0$. The computation of the coefficients $\alpha_{r+k}$ can be done through a careful application of the inverse bilateral Laplace transform.

    Alternatively, the second way to prove the claim for general forms is through an unilluminating computation using the definition.
    \begin{equation*}
      \begin{split}
        \mbf d(\mbf d\alpha)&=d(d\alpha)-s^{-1}ds\wedge\mf L_Xd\alpha-d\left(s^{-1}ds\wedge\mf L_X\alpha\right)+s^{-1}ds\wedge\mf L_X\left(s^{-1}ds\wedge\mf L_X\alpha\right)\\
        &=-s^{-1}ds\wedge\mf L_Xd\alpha-d(s^{-1})\wedge ds\wedge\mf L_X\alpha+s^{-1}ds\wedge d\mf L_X\alpha\\
        &\qquad+s^{-1}ds\wedge\mf L_X\left(s^{-1}ds\right)\wedge\mf L_X\alpha+s^{-1}ds\wedge s^{-1}ds\wedge\mf L_X\mf L_X\alpha\\
        &=-s^{-1}ds\wedge\mf L_Xd\alpha+s^{-1}ds\wedge d\mf L_X\alpha=0.
      \end{split}
    \end{equation*}
    Here, we used the facts that $d(d\alpha)=0$, $d$ satisfies a signed Leibniz rule, $\mf L_X$ satisfies the ordinary Leibniz rule, $ds\wedge ds=0$, $s^{-1}ds\in\Lambda^1_{-1}(\con)$ so $\mf L_X\left(s^{-1}ds\right)=0$, and $\mf L_X$ commutes with $d$.
  \end{proof}
\end{prop}

\begin{prop}\label{dantideriv}
  The operator $\mbf d$ is an anti-derivation. That is, if $\alpha\in\Lambda^k(\con)$ and $\beta\in\Lambda^{k'}(\con)$, then
  \begin{equation*}
    \mbf d(\alpha\wedge\beta)=\mbf d\alpha\wedge\beta+(-1)^k\alpha\wedge\mbf d\beta.
  \end{equation*}
  \begin{proof}[Proof for homogeneous forms]
    As before, the most illuminating computation is for homogeneous forms. Let $\alpha\in\Lambda^k_r(\con)$ and let $\beta\in\Lambda^{k'}_{r'}(\con)$. Let $\alpha=\gamma\,s^{r+k}$ and $\beta=\delta\,s^{r'+k'}$, so $\gamma$ and $\delta$ have total degree zero. Then $\alpha\wedge\beta=\gamma\wedge\delta\,s^{r+r'+k+k'}$, and $\gamma\wedge\delta$ has total degree zero. Thus, using the fact that $d$ is an anti-derivation, we find that
    \begin{equation*}
      \begin{split}
        \mbf d(\alpha\wedge\beta)&=d(\gamma\wedge\delta)s^{r+r'+k+k'}\\
        &=\left(d\gamma\wedge\delta+(-1)^k\gamma\wedge d\delta\right)\,s^{r+k}s^{r'+k'}\\
        &=d\gamma\,s^{r+k}\wedge\delta\,s^{r'+k'}+(-1)^k\gamma\,s^{r+k}\wedge d\delta\,s^{r'+k'}\\
        &=\mbf d\alpha\wedge\beta+(-1)^k\alpha\wedge\mbf d\beta.\qedhere
      \end{split}
    \end{equation*}
  \end{proof}
  \begin{proof}[Proof for general forms]
    For general forms that are not necessarily homogeneous, we can apply the same reasoning as in Proposition \ref{ddzero} to decompose a general form into its homogeneous components, or we can compute from the definition that
    \begin{equation*}
      \begin{split}
        \mbf d(\alpha\wedge\beta)&=d(\alpha\wedge\beta)-s^{-1}ds\wedge\mf L_X(\alpha\wedge\beta)\\
        &=d\alpha\wedge\beta+(-1)^k\alpha\wedge d\beta-s^{-1}ds\wedge\left(\mf L_X\alpha\wedge\beta+\alpha\wedge\mf L_X\beta\right)\\
        &=d\alpha\wedge\beta+(-1)^k\alpha\wedge d\beta\\
        &\qquad{}-s^{-1}ds\wedge\mf L_X\alpha\wedge\beta-(-1)^k\alpha\wedge s^{-1}ds\wedge\mf L_X\beta\\
        &=\left(d\alpha-s^{-1}ds\wedge\mf L_X\alpha\right)\wedge\beta+(-1)^k\alpha\wedge\left(d\beta-s^{-1}ds\wedge\mf L_X\beta\right)\\
        &=\mbf d\alpha\wedge\beta+(-1)^k\alpha\wedge\mbf d\beta,
      \end{split}
    \end{equation*}  
    as desired.
  \end{proof}
\end{prop}
\begin{cor}
  The operator $\mbf d$ commutes with multiplication by any power of $s$.
  \begin{proof}
    Since $s^r\in\Lambda^0_r(\con)$, we have that $\mbf d(s^r)=s^rd(1)=0$ by Proposition \ref{dhomog}. Then, by Proposition \ref{dantideriv},
    \begin{equation*}
      \mbf d(s^r\alpha)=\mbf d(s^r)\alpha+s^r\mbf d\alpha=s^r\mbf d\alpha.\qedhere
    \end{equation*}
  \end{proof}
\end{cor}
\begin{cor}\label{ddscommute}
  The operator $\mbf d$ commutes with the right wedge operation $\wedge ds$.
  \begin{proof}
    Since $ds\in\Lambda^1_0(\con)$, we have that
    \begin{equation*}
      \mbf d(ds)=d(ds)-s^{-1}ds\wedge\mf L_X(ds)=-s^{-1}ds\wedge ds=0.
    \end{equation*}
    By Lemma \ref{dantideriv}, for any $\alpha\in\Lambda^k(\con)$,
    \begin{equation*}
      \mbf d(\alpha\wedge ds)=\mbf d\alpha\wedge ds+(-1)^k\alpha\wedge\mbf d(ds)=\mbf d\alpha\wedge ds,
    \end{equation*}
    as desired.
  \end{proof}
\end{cor}
We conclude that, like $d$, $\mbf d$ sends vertical forms to vertical forms. In fact, as we will see in the proof of Theorem \ref{dcommutesvert}, when acting on vertical forms, $\mbf d$ is equal to $d$.
\begin{cor}\label{dvert}
  If $\alpha\in\Lambda^k(\con)^\perp$, then $\mbf d\alpha\in\Lambda^{k+1}(\con)^\perp$.
  \begin{proof}
    The space $\Lambda^k(\con)^\perp$ is defined to be $ds\wedge\Lambda^{k-1}(\con)$, which is equal to $\Lambda^{k-1}(\con)\wedge ds$, because commuting the order simply multiplies the form by $(-1)^{k-1}$. For a general $\alpha\wedge ds\in\Lambda^k(\con)^\perp$, by Corollary \ref{ddscommute}, we have
    \begin{equation*}
      \mbf d(\alpha\wedge ds)=\mbf d\alpha\wedge ds\in\Lambda^{k+1}(\con)^\perp,
    \end{equation*}
    as desired.
  \end{proof}
\end{cor}

Unlike $d$, the operator $\mbf d$ is also compatible with the horizontal subspace of forms.
\begin{prop}\label{djxcommute}
  The operator $\mbf d$ commutes with the right contraction $j_X$.
  \begin{proof}[Proof for homogenenous forms]
    Let $\alpha\in\Lambda^k_r(\con)$, and write $\alpha=\beta\,s^{r+k}$, so $\beta$ has total degree zero. On $k$-forms, $j_X=(-1)^{k-1}i_X$, so
    \begin{multline*}
      j_Xd\beta-dj_X\beta=(-1)^ki_Xd\beta-(-1)^{k-1}di_X\beta\\
      =(-1)^k(i_Xd\beta+di_X\beta)=(-1)^k\mf L_X\beta=0.
    \end{multline*}
    Since $\beta$ has total degree zero, $\mbf d\beta=d\beta$. Since $j_X$ lowers form degree and raises homogeneous degree, $j_X\beta$ also has total degree zero, so $\mbf dj_X\beta=dj_X\beta$. We conclude that
    \begin{equation*}
      j_X\mbf d\beta=\mbf dj_X\beta.
    \end{equation*}
    Multiplying the above equation by $s^{r+k}$ and using the fact that both $\mbf d$ and $j_X$ commute with multiplication by $s^{r+k}$, we find that
    \begin{equation*}
      j_X\mbf d\alpha=\mbf dj_X\alpha,
    \end{equation*}
    as desired.
  \end{proof}
  \begin{proof}[Proof for general forms.]
    As before, our focus is on homogeneous forms. However, if needed, we can decompose a general form into homogeneous components and apply the result for homogeneous forms. Alternatively, working from the definition and using the equation $j_Xd-dj_X=(-1)^k\mf L_X$ above, we find that, for $\alpha\in\Lambda^k(\con)$,
    \begin{equation*}
      \begin{split}
        &j_X\mbf d\alpha-\mbf dj_X\alpha\\
        &=j_Xd\alpha-dj_X\alpha-j_X\left(s^{-1}ds\wedge\mf L_X\alpha\right)+s^{-1}ds\wedge\mf L_Xj_X\alpha\\
        &=(-1)^k\mf L_X\alpha-(-1)^ki_X\left(s^{-1}ds\wedge\mf L_X\alpha\right)+(-1)^{k-1}s^{-1}ds\wedge\mf L_Xi_X\alpha\\
        &=(-1)^k\left(\mf L_X\alpha-i_X\left(s^{-1}ds\right)\mf L_X\alpha+s^{-1}ds\wedge i_X\mf L_X\alpha-s^{-1}ds\wedge\mf L_Xi_X\alpha\right)\\
        &=0,
      \end{split}
    \end{equation*}
    because $i_X\left(s^{-1}ds\right)=1$ and $\mf L_X\circ i_X=i_X\circ d\circ i_X=i_X\circ \mf L_X$ using the formula $\mf L_X=i_X\circ d+d\circ i_X$ and the fact that $i_X\circ i_X=0$.
  \end{proof}
\end{prop}
\begin{cor}\label{dhor}
  If $\alpha\in\Lambda^k(\con)^\top$, then $\mbf d\alpha\in\Lambda^{k+1}(\con)^\top$.
  \begin{proof}
    The space $\Lambda^k(\con)^\top$ is defined to be $i_X\Lambda^{k+1}(\con)$, which is equal to $j_X\Lambda^{k+1}(\con)$, because $j_X=(-1)^ki_X$. For a general $j_X\alpha\in\Lambda^k(\con)^\top$, we have
    \begin{equation*}
      \mbf dj_X\alpha=j_X\mbf d\alpha\in\Lambda^{k+1}(\con)^\top,
    \end{equation*}
    as desired.
  \end{proof}
\end{cor}

We finish this section by discussing the correspondence between the operator $\mbf d$ on the $\mbf P_r$ and $\mbf P_r^-$ spaces and the operator $d$ on the $\mc P_r$ and $\mc P_r^-$ spaces via the isomorphisms in Theorems \ref{prmthm} and \ref{prthm}.

\begin{prop}\label{drestrict}
  Let $i\colon T\to\con$ denote the inclusion of the simplex into the orthant, so the pullback $i^*$ denotes the restriction of forms. For $\alpha\in\Lambda^k(\con)$, we have
  \begin{equation*}
    i^*\mbf d\alpha=i^*d\alpha=di^*\alpha.
  \end{equation*}
  \begin{proof}
    The second equality is the naturality of the exterior derivative under pullback. The first equality follows from the definition, as follows.
    \begin{equation*}
      i^*\mbf d\alpha=i^*\left(d\alpha-s^{-1}ds\wedge\mf L_X\alpha\right)=i^*d\alpha-i^*(s^{-1})i^*(ds)\wedge i^*\mf L_X\alpha=i^*d\alpha,
    \end{equation*}
    because $i^*(ds)=0$.
  \end{proof}
\end{prop}

\begin{thm}\label{dcommutesvert}
  The operator $\mbf d$ maps $\mbf P_r\Lambda^k(\cO)$ into $\mbf P_{r-1}\Lambda^{k+1}(\cO)$, and the following diagram commutes.
  \begin{equation*}
    \begin{tikzcd}[column sep=large]
       \mbf P_r\Lambda^k(\cO) \arrow[d, "\mbf d"] \arrow[r, shift left, "i^*\circ j_X"] & \mc P_r\Lambda^k(T) \arrow[d, "d"] \arrow[l, shift left, "(\wedge ds)\circ h_r"]\\
       \mbf P_{r-1}\Lambda^{k+1}(\cO) \arrow[r, shift left, "i^*\circ j_X"]&\mc P_{r-1}\Lambda^{k+1}(T) \arrow[l, shift left, "(\wedge ds)\circ h_{r-1}"]
    \end{tikzcd}
  \end{equation*}
  \begin{proof}
    Let $\alpha\in\mbf P_r\Lambda^k(\cO)=\{\text{polynomial differential forms}\}\cap\Lambda^{k+1}_r(\con)^\perp$. Corollary \ref{dvert} tells us that $\mbf d$ sends vertical forms to vertical forms, and the fact that $\mbf d$ decreases homogeneous degree by one is clear from the definition of $\mbf d$. What is not apparent from the definition is that $\mbf d\alpha$ is a polynomial differential form, since the definition of $\mbf d$ involves dividing by $s$. But, since $\alpha$ is a vertical form, we have $ds\wedge\alpha=0$, and so
    \begin{equation*}
      \mbf d\alpha=d\alpha-s^{-1}ds\wedge\mf L_X\alpha=d\alpha-s^{-1}(r+k)ds\wedge\alpha=d\alpha.
    \end{equation*}
    Thus, indeed, if $\alpha$ is a polynomial differential form, then so is $\mbf d\alpha$.

    We now check that the diagram commutes. Indeed, by Proposition \ref{djxcommute}, $\mbf d$ commutes with $j_X$, so, applying Proposition \ref{drestrict}, we have that
    \begin{equation*}
      i^*j_X(\mbf d\alpha)=i^*(\mbf dj_X\alpha)=d(i^*j_X\alpha),
    \end{equation*}
    as desired.
  \end{proof}
\end{thm}

For the $\mbf P_r^-$ and $\mc P_r^-$ spaces, we must modify the above claim, as one would expect from the fact that $d$ fails to map $\mc P_r^-\Lambda^k(T)$ into $\mc P_{r-1}^-\Lambda^{k+1}(T)$, instead only mapping it into the larger space $\mc P_r^-\Lambda^{k+1}(T)$.

\begin{thm}
  The operator $s\mbf d$ maps $\mbf P_r^-\Lambda^k(\cO)$ into $\mbf P_r^-\Lambda^{k+1}(\cO)$, and the following diagram commutes.
  \begin{equation*}
    \begin{tikzcd}
       \mbf P_r^-\Lambda^k(\cO) \arrow[d, "s\mbf d"] \arrow[r, shift left, "i^*"] & \mc P_r^-\Lambda^k(T) \arrow[d, "d"] \arrow[l, shift left, "h_r"]\\
       \mbf P_r^-\Lambda^{k+1}(\cO) \arrow[r, shift left, "i^*"]&\mc P_r^-\Lambda^{k+1}(T) \arrow[l, shift left, "h_r"]
    \end{tikzcd}
  \end{equation*}
  \begin{proof}
    Let $\alpha\in\mbf P_r^-\Lambda^k(\cO)$. Recall that $\mbf P_r^-\Lambda^k(\cO)$ is the space of forms in $\Lambda^k_r(\con)^\top$ with polynomial coefficients. Similarly to earlier, Corollary \ref{dhor} tells us that $\mbf d$ sends horizontal forms to horizontal forms. Likewise, multiplication by $s$ sends horizontal forms to horizontal forms, since multiplication by a scalar field commutes with $i_X$. As before, $\mbf d$ decreases homogeneous degree by one, and multiplication by $s$ increases homogeneous degree by one. Thus, $s\,\mbf d\alpha\in\Lambda^{k+1}_r(\con)^\top$. To check that $s\,\mbf d\alpha$ is a polynomial differential form, we observe that, by definition,
    \begin{equation*}
      s\,\mbf d\alpha=s\,d\alpha-ds\wedge\mf L_X\alpha=s\,d\alpha-ds\wedge(r+k)\alpha,
    \end{equation*}
    which is clearly a polynomial differential form if $\alpha$ is. Thus, $s\,\mbf d\alpha\in\mbf P_r^-\Lambda^{k+1}(\cO)$, as desired.

    The diagram commutes by Proposition \ref{drestrict}. Indeed,
    \begin{equation*}
      i^*(s\,\mbf d\alpha)=(i^*s)(i^*\mbf d\alpha)=di^*\alpha,
    \end{equation*}
    since $i^*s=1$.
  \end{proof}
\end{thm}

\section{Integration and duality}\label{integrationsec}
In Arnold, Falk, and Winther \cite{afw06,afw10}, a key result leading to their geometric decomposition of the dual polynomial finite element spaces is that $\mc P_r^-\Lambda^k(T)$ is dual to $\mr{\mc P}_{r+k}\Lambda^{n-k}$ and $\mc P_r\Lambda^k(T)$ is dual to $\mr{\mc P}_{r+k+1}^-\Lambda^{n-k}(T)$ via the pairing $(a,b)\mapsto\int_Ta\wedge b$, where the ring over the $\mc P$ denotes those forms with vanishing trace when restricted to $\partial T$. We show that, in the language of $\mbf P_r$ and $\mbf P_r^-$ spaces, this duality occurs naturally via the Hodge star operator with respect to the nonstandard metric $g$.

In this section, we assume that $r\ge0$, though some of the results hold in greater generality.

\subsection{Integration}


\begin{dfn}
  Let $\mbf T\subset\cO$ denote the $(n+1)$-dimensional simplex
  \begin{equation*}
    \mbf T=\{x\in\cO\mid s\le1\}.
  \end{equation*}
\end{dfn}
The boundary of $\mbf T$ consists of the $n$-dimensional simplex $T=\{s=1\}$ along with the $n+1$ hyperplanes $x_i=0$ intersected with $\mbf T$.

\begin{prop}\label{integrationprop}
  Let $\mu\in\Lambda^{n+1}_r(\cO)$. Then
  \begin{equation*}
    (n+r+1)\int_{\mbf T}\mu=\int_Ti_X\mu.
  \end{equation*}
  \begin{proof}
    By Lemma \ref{ixlemma}, since $\mu$ is a top-level form,
    \begin{equation*}
      \int_{\mbf T}(n+r+1)\mu=\int_{\mbf T}\left(di_X\mu+i_X(0)\right)=\int_{\partial\mbf T}i_X\mu.
    \end{equation*}
    Consider the hyperplane $\mg_i=\{x_i=0\}$. On $\mg_i$, the outward-pointing vector field $X$ is tangent to $\mg_i$. As a consequence, $i_X\mu=0$ when restricted to $\mg_i$. Indeed, if $X_1,\dotsc,X_n$ are tangent to $\mg_i$, then
    \begin{equation*}
      i_X\mu(X_1,\dotsc,X_n)=\mu(X,X_1,\dotsc,X_n).
    \end{equation*}
    The $n+1$ vectors $X,X_1,\dotsc,X_n$ are all tangent to the $n$-dimensional space $\mg_i$ and are thus linearly dependent. Consequently, the above expression is zero. This leaves only one remaining boundary component of $\mbf T$, namely $T$, so we have
    \begin{equation*}
      \int_{\mbf T}(n+r+1)\mu=\int_{\partial\mbf T}i_X\mu=\int_{T}i_X\mu,
    \end{equation*}
    as desired.
  \end{proof}
\end{prop}
We now show that integrating forms over $T$ is, up to a nonzero constant, the same as integrating forms over $\mbf T$.
\begin{thm}
  Let $\alpha\in\Lambda^k_r(\cO)^\top$ and $\beta\in\Lambda^{n+1-k}_{r'}(\cO)^\perp$. Let $a\in\Lambda^k(T)$ and $b\in\Lambda^{n-k}(T)$ be the forms corresponding to $\alpha$ and $\beta$ via the isomorphisms in Corollaries \ref{homtcor} and \ref{vertcor}. Then
  \begin{equation*}
    (-1)^n(n+r+r'+1)\int_{\mbf T}\alpha\wedge\beta=\int_Ta\wedge b.
  \end{equation*}
  \begin{proof}
    Since $\alpha\wedge\beta\in\Lambda^{n+1}_{r+r'}(\cO)$, Proposition \ref{integrationprop},
    \begin{equation*}
      (-1)^n(n+r+r'+1)\int_{\mbf T}\alpha\wedge\beta=(-1)^n\int_Ti_X(\alpha\wedge\beta)=\int_Tj_X(\alpha\wedge\beta).
    \end{equation*}
    Because $\alpha$ is a horizontal form, $j_X\alpha=0$. Thus, by the appropriate signed Leibniz rule for $j_X$, we have
    \begin{equation*}
      \int_Tj_X(\alpha\wedge\beta)=\int_T\alpha\wedge j_X\beta=\int_Ta\wedge b,
    \end{equation*}
    because $a=i^*\alpha$ and $b=i^*(j_X\alpha)$.
  \end{proof}
\end{thm}

\subsection{The Hodge star with respect to $g$}
\begin{dfn}
  Let the Hodge star with respect to $g$ and the Euclidean volume, denoted $*_g$, be defined pointwise by
  \begin{equation*}
    \alpha\wedge *_g\beta=\langle\alpha,\beta\rangle_g\vol,
  \end{equation*}
  where $\vol$ denotes the standard Euclidean volume form (as opposed to the volume form induced by $g$), and $\alpha,\beta\in\Lambda^k(\cO)$.
\end{dfn}

Note that $*_g$ is \emph{not} the standard Hodge star operation on $(\mbf O,g)$ viewed as a Riemannian manifold, because the standard definition would require that we use the volume form $\mu_g=\sqrt{\det g}\,\vol$ corresponding to the metric $g$, not the Euclidean volume form $\vol$. However, the relationship between the two is simple; see the proof of Proposition \ref{starstarxyz}.

One may reasonably ask if the degeneracy of the inner product on the boundary of $\cO$ leads to any issues with this definition. It does not. Indeed, given $\beta$, we view $\alpha\mapsto\langle\alpha,\beta\rangle_g\vol$ as a linear functional $\Lambda^k_x(\cO)\to\Lambda^{n+1}_x(\cO)\cong\mb R$. A standard property of the exterior algebra is that the wedge pairing $\Lambda^k_x(\cO)\times\Lambda^{n+1-k}_x(\cO)\xrightarrow\wedge\Lambda^{n+1}_x(\cO)$ is a perfect pairing, so the linear functional $\alpha\mapsto\langle\alpha,\beta\rangle_g\vol$ can be represented by a map of the form $\alpha\mapsto\alpha\wedge\gamma$ for a unique $\gamma\in\Lambda^{n+1-k}_x(\cO)$. We let $*_g\beta$ denote this form $\gamma$.

\begin{eg}
  If $n=2$, then one can evaluate $*_g$ as follows.
  \begin{equation*}
    \begin{array}{c|c}
      \beta&*_g\beta\\\hline
      1&dx\wedge dy\wedge dz\\
      dx&x\,dy\wedge dz\\
      dx\wedge dy&xy\,dz\\
      dx\wedge dy\wedge dz&xyz
    \end{array}
  \end{equation*}
\end{eg}

The degeneracy of $g$ does lead to the unusual behavior that, at a point $x$ on the boundary of $\cO$, it may happen that $*_g\beta_x=0$ even though $\beta_x\neq0$. This behavior is by design, as shown in the following proposition.
\begin{dfn}
  Let $\mg_i=\{x\in\cO\mid x_i=0\}$ be the $i$th hyperplane comprising the boundary of $\cO$, and let $k_i\colon\mg_i\hookrightarrow\cO$ be the inclusion. 
\end{dfn}
\begin{prop}\label{starnotrace}
  If $\beta\in\Lambda^k(\cO)$, then $*_g\beta\in\mr\Lambda^{n+1-k}(\cO)$.
  \begin{proof}
    It suffices to show that $k_i^*(*_g\beta)=0$ for all $i$.

    Let $x\in\mg_i$. Let $\gamma\in\Lambda^{k-1}_x$. Note that, at the point $x$, $\langle dx_i,\phi\rangle_g=0$ for any $\phi\in T^*_x\mb R^{n+1}$. That is, $dx_i$ is in the kernel of $g$. Thus, by the definition of $*_g$, we have
    \begin{equation*}
      \gamma\wedge dx_i\wedge *_g\beta_x=\langle\gamma\wedge dx_i,\beta_x\rangle_g\vol=0.
    \end{equation*}
    Here we are using the fact that $dx_i$ is in the kernel of $g$, so when $g$ is extended to the exterior algebra, the form $\gamma\wedge dx_i$ must also be in the kernel of $g$.

    Because $\gamma\wedge dx_i\wedge *_g\beta_x=0$ for all choices of $\gamma\in\Lambda^{k-1}_x$, we conclude by properties of the exterior algebra that $dx_i\wedge *_g\beta_x=0$. Using the formula $(dx_i\wedge)\circ i_{\pp{x_i}}+i_{\pp{x_i}}\circ dx_i=1$, we conclude that $*_g\beta_x=dx_i\wedge i_{\pp{x_i}}(*_g\beta_x)$. Since $k_i^*(dx_i)=0$, we conclude that $k_i^*(*_g\beta_x)=0$, as desired.
  \end{proof}
\end{prop}
Note that, depending on regularity assumptions, the converse may not be true. For example $\sqrt x\,dy\wedge dz$ has vanishing trace on $\partial\cO$, but $*_g^{-1}(\sqrt x\,dy\wedge dz)=\frac1{\sqrt x}\,dx$, which is defined on $\mbf O$ but is discontinuous at the $x=0$ boundary. As we will show in Proposition \ref{hodgepolyiso}, the converse does hold for polynomials.

We can refine Proposition \ref{starnotrace} as follows.
\begin{prop}\label{polystarnotrace}
  If $\beta\in\Lambda^k_r(\cO)$, then $*_g\beta\in\mr\Lambda^{n+1-k}_{r+k}(\cO)$. If $\beta\in\mbf H_r\Lambda^k(\cO)$, then $*_g\beta\in\mr{\mbf H}_{r+k}\Lambda^{n+1-k}(\cO)$.
  \begin{proof}
    Given what we have already proved, the content of the claim is that $*_g$ raises homogeneous degree by $k$, and that if $\beta$ is a polynomial differential form, then so is $*_g\beta$.

    We have been using $g$ to denote the inner product on differential forms of any degree. For this proof, however, we will need to be more specific. Let $g_k$ denote $g$ on forms of degree $k$. From the definition of $g$, we have
    \begin{equation*}
      g_1=\sum_{i=1}^{n+1}x_i\pp{x_i}\otimes\pp{x_i}.
    \end{equation*}
    Note that the coefficients of $g_1$ are polynomials of degree one. Viewing $\bigwedge^kT^*\cO$ as the alternating subset of $\bigotimes^kT^*\cO$ with an appropriate normalization, we can write $g_k$ as the $k$-fold tensor product
    \begin{equation*}
      g_k=g_1\otimes\dotsb\otimes g_1.
    \end{equation*}
    The coefficients of $g_k$ are thus polynomials of degree $k$.

    We now consider the formula $\alpha\wedge *_g\beta=\langle\alpha,\beta\rangle_g\vol$. Since the Hodge star is a pointwise (tensorial) operation, we can assume without loss of generality that $\alpha$ is a constant $k$-form. Then, $\langle\alpha,\beta\rangle_g$ is a sum of products of the coefficients of $\alpha$, $\beta$, and $g_k$. Thus, $\langle\alpha,\beta\rangle_g$ is homogeneous of degree $r+k$ and is a polynomial if the coefficients of $\beta$ are polynomials.

    Thus, $\alpha\wedge*_g\beta$ is a form of homogeneous degree $r+k$ and a polynomial differential form if $\beta$ is. Since $\alpha$ is a constant $k$-form, $*_g\beta$ is thus forced to be homogeneous of degree $r+k$ and a polynomial differential form if $\beta$ is.
  \end{proof}
\end{prop}

\begin{prop}\label{starstarvanish}
  If $\beta\in\mr\Lambda^k(\cO)$, then $*_g\beta_x=0$ for all $x\in\partial\cO$. In other words, all components of $*_g\beta$ vanish on the boundary, not just the tangential ones.
  \begin{proof}
    Let $x\in\partial\cO$. For clarity of exposition, we first consider the case where $x$ is on exactly one boundary face $\mg_i=\{x_i=0\}$. At this point $x$, the condition that $\beta\in\mr\Lambda^k(\cO)$ can be written as $k_i^*\beta_x=0$. Either using a basis or with an argument analogous to Proposition \ref{htiso} except with $\mg_i$ instead of $T$, we can show that this condition implies that $\beta_x=dx_i\wedge\gamma$ for some $\gamma\in\Lambda^{k-1}_x(\cO)$. As before, the fact that $dx_i$ is in the kernel of $g$ at $x$ implies that
    \begin{equation*}
      \langle\alpha,dx_i\wedge\gamma\rangle_g=0
    \end{equation*}
    for all $\alpha\in\Lambda^k_x(\cO)$. Hence,
    \begin{equation*}
      \alpha\wedge*_g\beta_x=\langle\alpha,dx_i\wedge\gamma\rangle_g\vol=0
    \end{equation*}
    for all $\alpha\in\Lambda^k_x(\cO)$, from which we can conclude that $*_g\beta_x=0$, as desired.

    For a general $x\in\partial\cO$ that may potentially be on multiple faces, we can obtain the result by continuity of $g$ and hence continuity of $*_g\beta$. Alternatively, assume that $x$ is on the intersection of several boundary faces $\Gamma_{i_1},\dotsc,\Gamma_{i_m}$. The condition that $\beta\in\mr\Lambda^k(\cO)$ is weaker in this situation, since it only requires that $\beta_x(X_1,\dotsc,X_k)=0$ for vectors $X_k$ that are tangent to \emph{all} the hyperplanes $\Gamma_{i_1},\dotsc,\Gamma_{i_m}$. Under this weaker hypothesis, we can still use our understanding of the exterior algebra to conclude that
    \begin{equation*}
      \beta_x=dx_{i_1}\wedge\gamma_1+\dotsb+dx_{i_m}\wedge\gamma_m
    \end{equation*}
    for some forms $\gamma_1,\dotsc,\gamma_m\in\Lambda^{k-1}_x(\cO)$. But, because $dx_{i_1}$ through $dx_{i_m}$ are \emph{all} in the kernel of $g$ at this point $x$, we still conclude that $\langle\alpha,\beta_x\rangle_g=0$ for all $\alpha$ and complete the proof as before.
  \end{proof}
\end{prop}

\begin{dfn}
  Let $p$ denote the product $x_1\dotsm x_{n+1}$.
\end{dfn}
The previous proposition is not too surprising given the following formula, along with the fact that $p$ vanishes on $\partial\cO$.
\begin{prop}\label{starstarxyz}
  On $k$-forms, we have the formula $*_g\circ *_g=(-1)^{k(n+1-k)}p$. That is, for any $\beta\in\Lambda^k(\cO)$, we have
  \begin{equation*}
    *_g(*_g\beta)=(-1)^{k(n+1-k)}p\beta.
  \end{equation*}
  \begin{proof}
    On the interior $\mbf O$, the metric $g$ is nondegenerate, and one can compute that the volume form corresponding to $g$ is $\mu_g=\frac1{\sqrt p}\,\vol$. Indeed,
    \begin{equation*}
      \langle\vol,\vol\rangle_g=\langle dx_1\wedge\dotsb\wedge dx_{n+1},dx_1\wedge\dotsb\wedge dx_{n+1}\rangle_g=x_1\dotsm x_{n+1}=p,
    \end{equation*}
    so $\langle\mu_g,\mu_g\rangle=\left\langle\frac1{\sqrt p}\vol,\frac1{\sqrt p}\vol\right\rangle=1$, as desired.

    Viewing $(\mbf O,g)$ as a Riemannian manifold, let $*_g'$ denote the standard Hodge star on Riemannian manifolds defined by $\alpha\wedge *_g'\beta=\langle\alpha,\beta\rangle_g\mu_g$. From standard references, this Hodge star operation satisfies $*_g'(*_g'\beta)=(-1)^{k(n+1-k)}\beta$. Compare this formula to our definition $\alpha\wedge *_g\beta=\langle\alpha,\beta\rangle_g\vol$. Because $\vol=\sqrt p\,\mu_g$, we have that $*_g=\sqrt p\,*_g'$. Hence $*_g\circ *_g=p(*_g'\circ *_g')$, and so we conclude that $*_g\circ *_g=p(-1)^{k(n+1-k)}$ on $\mbf O$.

    To extend this result to $\cO$, we can either use continuity of $g$ and hence $*_g$, or we can apply Propositions \ref{starnotrace} and \ref{starstarvanish} to conclude that $*_g\circ *_g=0$ on $\partial\cO$ and noting that $p=0$ on $\partial\cO$.
  \end{proof}
\end{prop}

\subsection{Duality}

We are now ready to show the converse of Proposition \ref{starnotrace} for polynomials. In other words, if $\alpha$ is a polynomial with vanishing restriction to the boundary, then $*_g^{-1}\alpha$ is a polynomial.
\begin{prop}\label{hodgepolyiso}
  The Hodge star with respect to $g$ induces an isomorphism between the following spaces of polynomial differential forms.
  \begin{equation*}
    \begin{tikzcd}
      \mbf H_r\Lambda^k(\cO)\arrow[r, shift left, "*_g"] &\mr{\mbf H}_{r+k}\Lambda^{n+1-k}(\cO). \arrow[l, shift left]
    \end{tikzcd}
  \end{equation*}
  \begin{proof}
    By Proposition \ref{polystarnotrace}, the map $*_g$ does indeed send the space $\mbf H_r\Lambda^k(\cO)$ to the space $\mr{\mbf H}_{r+k}\Lambda^{n+1-k}(\cO)$. To construct an inverse, the key idea is that $*_g^{-1}=\pm\frac1p\,*_g$ by Proposition \ref{starstarxyz}. However, there is no reason to believe a priori that dividing by $p$ yields a polynomial.

    Let $\alpha\in\mr{\mbf H}_{r+k}\Lambda^{n+1-k}(\cO)$. Then $*_g\alpha\in\mbf H_{r+n+1}\Lambda^k(\cO)$, and, moreover, by Proposition \ref{starstarvanish}, $*_g\alpha_x=0$ for all $x\in\partial\cO$. In other words, all components of $*_g\alpha$ vanish on the boundary of $\cO$, not just the tangential ones as with restriction. Consequently, all of the polynomial coefficients of $*_g\alpha$ vanish on $\partial\cO$. That is, these polynomial coefficients vanish when $x_i=0$ for any $i$. A polynomial that vanishes when $x_i=0$ must be divisible by $x_i$. Thus, the coefficients of $*_g\alpha$ are divisible by $p=x_1\dotsm x_{n+1}$ as polynomials. Hence, $*_g\alpha=p\beta$ for some $\beta\in\mbf H_r\Lambda^k(\cO)$. Applying $*_g$ to both sides, we find that $(-1)^{k(n+1-k)}p\alpha=*_g*_g\alpha=p*_g\beta$, so $\alpha=*_g(-1)^{k(n+1-k)}\beta$ for $\beta\in\mbf H_r\Lambda^k(\cO)$, as desired.
  \end{proof}
\end{prop}

As discussed previously, the statement above fails to hold without the polynomial assumption. For example, with $n=2$, the form $\alpha=\sqrt x\,dy\wedge dz$ has zero tangential trace on $\partial\cO$, but $*^{-1}_g\alpha=\frac1p*_g\alpha=\frac1{\sqrt x}\,dx$, which is defined on $\mbf O$ but does not extend continuouly to the $x=0$ boundary of $\cO$. However, the statement does hold for smooth differential forms; the key step that $*_g\alpha$ is divsible by $p$ follows from Taylor's theorem. Weaker assumptions may suffice.

We now have the tools to show that the integration pairing induces a duality relationship between these polynomial spaces.

\begin{thm}\label{hpairing}
  The pairing
  \begin{equation*}
    (\alpha,\beta)\mapsto\int_{\mbf T}\alpha\wedge\beta
  \end{equation*}
  is nondegenerate for $(\alpha,\beta)\in\mbf H_r\Lambda^k(\cO)\times\mr{\mbf H}_{r+k}\Lambda^{n+1-k}(\cO)$.

  That is, if $\alpha\in\mbf H_r\Lambda^k(\cO)$ and $\int_{\mbf T}\alpha\wedge\beta$ is zero for all $\beta\in\mr{\mbf H}_{r+k}\Lambda^{n+1-k}(\cO)$, then $\alpha$ is zero, and likewise if $\beta\in\mr{\mbf H}_{r+k}\Lambda^{n+1-k}(\cO)$ and $\int_{\mbf T}\alpha\wedge\beta$ is zero for all $\alpha\in\mbf H_r\Lambda^k(\cO)$, then $\beta$ is zero.

  As a consequence, this pairing induces the isomorphism
  \begin{align*}
    \left(\mbf H_r\Lambda^k(\cO)\right)^*&\cong\mr{\mbf H}_{r+k}\Lambda^{n+1-k}(\cO).
  \end{align*}

  \begin{proof}
    Let $\alpha\in\mbf H_r\Lambda^k(\cO)$, and assume that $\int_{\mbf T}\alpha\wedge\beta=0$ for all $\beta\in\mr{\mbf H}_{r+k}\Lambda^{n+1-k}$. By Proposition \ref{hodgepolyiso}, $*_g\alpha\in\mr{\mbf H}_{r+k}\Lambda^{n+1-k}$. Thus,
    \begin{equation}\label{aazeroeqn}
      0=\int_{\mbf T}\alpha\wedge *_g\alpha=\int_{\mbf T}\langle\alpha,\alpha\rangle_g\vol.
    \end{equation}
    We conclude that $\alpha=0$, though we must be slightly careful because of the degeneracy of $g$. We have that $\langle\alpha,\alpha\rangle_g\vol\ge0$ even on the boundary where $g$ is degenerate, so, from the fact that the integral is zero, we can conclude that $\langle\alpha,\alpha\rangle_g=0$ almost everywhere on $\mbf T$. On the interior of $\mbf T$, the inner product $g$ is non-degenerate, so we conclude that $\alpha=0$ almost everywhere on the interior of $\mbf T$. Since $\alpha$ is continuous, we conclude that $\alpha=0$ on all of $\mbf T$.

    Now let $\beta\in\mr{\mbf H}_{r+k}\Lambda^{n+1-k}(\cO)$, and assume that $\int_{\mbf T}\alpha\wedge\beta=0$ for all $\alpha\in\mbf H_r\Lambda^k(\cO)$. By Proposition \ref{hodgepolyiso}, $\beta=*_g\alpha$ for some $\alpha\in\mbf H_r\Lambda^k(\cO)$. Thus, equation \eqref{aazeroeqn} holds, and we conclude that $\alpha=0$, so $\beta=0$.

    The duality relationship
    \begin{align*}
      \left(\mbf H_r\Lambda^k(\cO)\right)^*&\cong\mr{\mbf H}_{r+k}\Lambda^{n+1-k}(\cO).
    \end{align*}
    follows from the fact that, in finite dimensions, a nondegenerate pairing is perfect.
  \end{proof}
\end{thm}


For the corresponding theorem for the $\mbf P_r$ and $\mbf P_r^-$ spaces, we need the following proposition.
\begin{prop}\label{starverthor}
  The Hodge star $*_g$ interchanges vertical and horizontal forms. That is, if $\beta\in\Lambda^k(\con)^\top$, then $*_g\beta\in\Lambda^{n+1-k}(\con)^\perp$, and if $\beta\in\Lambda^k(\con)^\perp$, then $*_g\beta\in\Lambda^{n+1-k}(\con)^\top$.
  \begin{proof}
    If $\beta\in\Lambda^k(\con)^\top$, then by Proposition \ref{dsixsplitting}, $\beta$ is pointwise $g$-orthogonal to $\Lambda^k(\con)^\perp$. Thus, for all $\gamma\in\Lambda^{k-1}(\con)$, we have $\langle\gamma\wedge ds,\beta\rangle_g=0$. Therefore,
    \begin{equation*}
      \gamma\wedge ds\wedge *_g\beta=\langle\gamma\wedge ds,\beta\rangle_g\vol=0
    \end{equation*}
    for all $\gamma\in\Lambda^{k-1}(\con)$. By properties of the exterior algebra, we conclude that $ds\wedge *_g\beta$ must itself be zero, so $*_g\beta\in\Lambda^{n+1-k}(\con)^\perp$, as desired.

    Conversely, if $\beta\in\Lambda^k(\con)^\perp$, then by Proposition \ref{dsixsplitting}, $\beta$ is pointwise $g$-orthogonal to $\Lambda^k(\con)^\top$. In other words, for all $\gamma\in\Lambda^{k+1}(\con)$, we have $\langle i_X\gamma,\beta\rangle_g=0$. Thus, for all $\gamma\in\Lambda^{k+1}(\con)$,
    \begin{equation*}
      0=\langle i_X\gamma,\beta\rangle_g\vol=i_X\gamma\wedge*_g\beta.
    \end{equation*}
    Note that $\gamma\wedge*_g\beta$ is an $(n+2)$-form in an $(n+1)$-dimensional space, and is hence zero. Thus, the signed Leibniz rule $i_X(\gamma\wedge*_g\beta)=i_X\gamma\wedge*_g\beta+(-1)^{k+1}\gamma\wedge i_X(*_g\beta)$ tells us that
    \begin{equation*}
      0=(-1)^ki_X\gamma\wedge*_g\beta=\gamma\wedge i_X(*_g\beta).
    \end{equation*}
    Since this exterior product vanishes for all $(k+1)$-forms $\gamma$, we conclude by properties of the exterior algebra that $i_X(*_g\beta)$ is itself zero, so $*_g\beta\in\Lambda^{n+1-k}(\con)^\top$, as desired.
  \end{proof}
\end{prop}
\begin{cor}\label{hodgevhiso}
  The Hodge star with respect to $g$ induces an isomorphism between the following spaces of forms on the interior of the orthant $\mbf O$.
  \begin{equation*}
    \begin{tikzcd}[row sep=0]
      \Lambda^k(\mbf O)^\top\arrow[r, shift left, "*_g"] &\Lambda^{n+1-k}(\mbf O)^\perp, \arrow[l, shift left]\\
      \Lambda^k(\mbf O)^\perp\arrow[r, shift left, "*_g"] &\Lambda^{n+1-k}(\mbf O)^\top. \arrow[l, shift left]\\
    \end{tikzcd}
  \end{equation*}
  \begin{proof}
    Proposition \ref{starverthor} shows that $*_g$ has the specified targets. For the inverse map, we use the fact that $*_g^{-1}=\pm\frac1p*_g$ from Proposition \ref{starstarxyz}. Since we are restricting ourselves to the interior $\mbf O$, the scalar function $\pm\frac1p$ is always defined. The decomposition into vertical and horizontal forms is defined pointwise, and is hence unaffected by multiplication by a scalar function. Thus, because $*_g$ sends horizontal forms to vertical forms and vertical forms to horizontal forms, so does $*_g^{-1}$.
  \end{proof}
\end{cor}
\begin{prop}\label{hodgepiso}
  Let $r\ge0$. The Hodge star with respect to $g$ induces an isomorphism between the following spaces of polynomial differential forms.
  \begin{equation*}
    \begin{tikzcd}[row sep=0,
      /tikz/column 1/.append style={anchor=base east}
      ,/tikz/column 2/.append style={anchor=base west}
      ]
      \mbf P_r^-\Lambda^k(\cO)\arrow[r, shift left, "*_g"] &\mr{\mbf P}_{r+k}\Lambda^{n-k}(\cO), \arrow[l, shift left]\\
      \mbf P_r\Lambda^k(\cO)\arrow[r, shift left, "*_g"] &\mr{\mbf P}_{r+k+1}^-\Lambda^{n-k}(\cO). \arrow[l, shift left]
    \end{tikzcd}
  \end{equation*}
  \begin{proof}
    A form that is vertical on $\mbf O$ and defined on $\con$ must be vertical on all of $\con$ by continuity, and likewise for horizontal forms. With this in mind, we can write
    \begin{align*}
      \mbf P_r^-\Lambda^k(\cO)&=\mbf H_r\Lambda^k(\cO)\cap\Lambda^k(\mbf O)^\top,\\
      \mbf P_r\Lambda^k(\cO)&=\mbf H_r\Lambda^{k+1}(\cO)\cap\Lambda^{k+1}(\mbf O)^\perp.
    \end{align*}
    In other words, $\mbf P_r^-\Lambda^k(\cO)$ contains those homogeneous polynomial $k$-forms of degree $r$ that are horizontal on $\mbf O$, and $\mbf P_r\Lambda^k(\cO)$ contains those homogeneous polynomial $(k+1)$-forms of degree $r$ that are vertical on $\mbf O$.

    With this understanding, the isomorphisms we seek to prove are
    \begin{equation*}
      \begin{tikzcd}[row sep=0,
        /tikz/column 1/.append style={anchor=base east}
        ,/tikz/column 2/.append style={anchor=base west}
        ]
        \mbf H_r\Lambda^k(\cO)\cap\Lambda^k(\mbf O)^\top\arrow[r, shift left, "*_g"] &\mr{\mbf H}_{r+k}\Lambda^{n+1-k}(\cO)\cap\Lambda^{n+1-k}(\mbf O)^\perp, \arrow[l, shift left]\\
        \mbf H_r\Lambda^{k+1}(\cO)\cap\Lambda^{k+1}(\mbf O)^\perp\arrow[r, shift left, "*_g"] &\mr{\mbf H}_{r+k+1}\Lambda^{n-k}(\cO)\cap\Lambda^{n-k}(\mbf O)^\top. \arrow[l, shift left]
      \end{tikzcd}
    \end{equation*}
    But these isomorphisms follow directly from the isomorphisms in Proposition \ref{hodgepolyiso} and Corollary \ref{hodgevhiso}.
  \end{proof}
\end{prop}
\begin{thm}\label{pairingthm}
  The pairing
  \begin{equation*}
    (\alpha,\beta)\mapsto\int_{\mbf T}\alpha\wedge\beta
  \end{equation*}
  is nondegenerate for either of
  \begin{align*}
    (\alpha,\beta)&\in \mbf P_r^-\Lambda^k(\cO)\times\mr{\mbf P}_{r+k}\Lambda^{n-k}(\cO),\\
    (\alpha,\beta)&\in\mbf P_r\Lambda^k(\cO)\times\mr{\mbf P}_{r+k+1}^-\Lambda^{n-k}(\cO).
  \end{align*}
  As a consequence, this pairing induces the isomorphisms
  \begin{align*}
    \left(\mbf P_r^-\Lambda^k(\cO)\right)^*&\cong\mr{\mbf P}_{r+k}\Lambda^{n-k}(\cO),&
    \left(\mbf P_r\Lambda^k(\cO)\right)^*&\cong\mr{\mbf P}_{r+k+1}^-\Lambda^{n-k}(\cO).
  \end{align*}
  \begin{proof}
    The proof is identical to the proof of Theorem \ref{hpairing}, except with Proposition \ref{hodgepiso} in place of Proposition \ref{hodgepolyiso}.
  \end{proof}
\end{thm}



\begin{cor}
  The pairing
  \begin{equation*}
    (a,b)\mapsto\int_Ta\wedge b
  \end{equation*}
  is nondegenerate for either of
  \begin{align*}
    (a,b)&\in \mc P_r^-\Lambda^k(T)\times\mr{\mc P}_{r+k}\Lambda^{n-k}(T),\\
    (a,b)&\in\mc P_r\Lambda^k(T)\times\mr{\mc P}_{r+k+1}^-\Lambda^{n-k}(T).
  \end{align*}
  \begin{proof}
    Theorems \ref{prmthm} and \ref{prthm} tell us that the $\mbf P_r^-$ and $\mbf P_r$ spaces are isomorphic to the $\mc P_r^-$ and $\mc P_r$ spaces, and by Proposition \ref{prnotrace}, these isomorphisms send the $\mr{\mbf P}_r^-$ and $\mr{\mbf P}_r$ spaces to the $\mr{\mc P}_r^-$ and $\mr{\mc P}_r$ spaces, respectively. Proposition \ref{integrationprop} tells us that, via these isomorphisms, the pairings $\int_{\mbf T}\alpha\wedge\beta$ and $\int_Ta\wedge b$ are related by a nonzero constant. Hence, if one pairing is nondegenerate, so is the other.
  \end{proof}  
\end{cor}

We thus obtain a new proof of Arnold, Falk, and Winther's key lemma \cite[Lemma 5.6]{afw10}.
\begin{cor}\label{afwcor}
  The pairing
  \begin{equation*}
    (a,b)\mapsto\int_Ta\wedge b
  \end{equation*}
  induces the isomorphisms
  \begin{align*}
    \left(\mr{\mc P}_r\Lambda^k(T)\right)^*&\cong\mc P_{r+k-n}^-\Lambda^{n-k}(T),&\left(\mr{\mc P}_r^-\Lambda^k(T)\right)^*&\cong\mc P_{r+k-n-1}\Lambda^{n-k}(T).
  \end{align*}
  \begin{proof}
    In finite dimensions, a nondegenerate pairing is perfect, so
    \begin{align*}
      \left(\mr{\mc P}_{r+k}\Lambda^{n-k}(T)\right)^*&\cong\mc P_r^-\Lambda^k(T),&\left(\mr{\mc P}_{r+k+1}^-\Lambda^{n-k}(T)\right)^*&\cong\mc P_r\Lambda^k(T).
    \end{align*}
    Reindexing by replacing $k$ with $n-k$ and $r$ with either $r+k-n$ or $r+k-n-1$, respectively, we obtain the desired statement.
  \end{proof}
\end{cor}

\section{Future work}\label{futuresec}
\begin{question}
  The $\mc P_r\Lambda^k$ and $\mc P^-_r\Lambda^k$ families give finite element methods for simplicial meshes. For parallelotope meshes, Arnold and Awanou construct the serendipity elements $\mc S_r\Lambda^k$ \cite{aa11}. Can we make an analogous construction of an $\mbf S_r\Lambda^k$ space for the serendipity elements, with a simple proof of the duality relationship corresponding to Theorem \ref{pairingthm} in this context? 
\end{question}

\begin{question}
  To understand how solutions given by finite element methods approximate the true solution, Falk and Winther construct local bounded cochain projections from the space of $L^2$ forms with $L^2$ exterior derivatives to the $\mc P$ and $\mc P^-$ spaces \cite{fw14}. Could the nonstandard inner product $g$ lead to another construction of a local bounded cochain projection?
\end{question}

\section{Acknowledgements} I would like to thank my postdoctoral mentor Ari Stern for the valuable feedback on this paper.

\bibliographystyle{plain}
\bibliography{fem}

\end{document}